\documentclass[11pt]{article}

\usepackage{amssymb,amsmath,amsthm,amsopn,amsfonts}
\usepackage[T1]{fontenc}
\usepackage{ae,aecompl}
\usepackage{graphicx}
\usepackage{bm}
\usepackage{fancybox,layout,color}
\usepackage{fancyhdr}

\usepackage{scrtime}

\usepackage{setspace}
\usepackage{url}
\usepackage{float}
\usepackage{placeins}

\numberwithin{equation}{section}
\definecolor{purple}{RGB}{160,32,40}

\usepackage{ulem}  % for the command \sout{}
\usepackage[makeroom]{cancel}  % for the command \varnothing{}

\newtheorem{teo}{Theorem}[section]
\newtheorem{nota}[teo]{Remark}

\newtheorem{coro}[teo]{Corollary}
\newtheorem{defi}[teo]{Definition}
\newtheorem{prop}[teo]{Proposition}
\newtheorem{lema}[teo]{Lemma}

\newcommand{\R}{\ensuremath{\mathbb{R}} }

\newcommand{\dist}{\mathrm{dist}}

\newcommand{\epi}{\mathrm{epi}}
\newcommand{\Aff}{\mathrm{Aff}}

\newcommand{\be}{\begin{eqnarray}}
\newcommand{\ee}{\end{eqnarray}}

\DeclareMathOperator{\co}{\mathsf{co}}

\title{ Compensated Convexity Methods\\ for Approximations and Interpolations\\
of Sampled Functions in Euclidean Spaces: \\ Theoretical Foundations}

\author{\normalsize  Kewei Zhang\thanks{School  of Mathematical Sciences,
University of Nottingham, University Park, Nottingham, NG7 2RD, UK},\,
Elaine Crooks\thanks{Department of Mathematics, Swansea University,
Singleton Park, Swansea, SA2 8PP, UK}
\,
and Antonio Orlando\thanks{CONICET, FACET,
Universidad Nacional de Tucum\'an, Argentina}
}

\date{ }

\begin{document}

%%%%%%%%%%%%%%%%%%%%%%%%%%%%%%%

\maketitle

%%%%%%%%%%%%%%%%%%%%%%%%%%%%%%%

\singlespacing
\pagestyle{fancy}
\fancyhead{}
% \fancyfoot[OR,ER]{\tiny \today\,\, \thistime}
\cfoot{\thepage}

%%%%%%%%%%%%%%%%%%%%%%%%%%%%%%%%

\begin{abstract}
We  introduce Lipschitz continuous and $C^{1,1}$ geometric approximation and interpolation 
methods for sampled bounded uniformly continuous functions over compact sets and over complements of bounded open sets
in $\mathbb{R}^n$ by using compensated convex transforms.  
Error estimates are provided for the approximations of bounded uniformly continuous functions, 
of Lipschitz functions, and 
of $C^{1,1}$ functions. We also prove that our approximation methods, which  are differentiation and integration free and not sensitive to sample type, are stable with respect
to  the Hausdorff distance between samples.
\end{abstract}

\medskip
\footnotesize{

\medskip
{\bf 2000 Mathematics Subjects Classification number}:
90C25, 90C26, 49J52, 52A41, 65K10

\medskip
{\bf Email}:  kewei.zhang@nottingham.ac.uk,    e.c.m.crooks@swansea.ac.uk, aorlando@herrera.unt.edu.ar

}

\normalsize

%%%%%%%%%%%%%%%%%%%%%%%%%%%%%%%%%%%%%%%%%%%%%%%%%%%%%%%%%%%%%%%%%%%%%%%%%%%%%%%%%%%%%%%%%%%%%%%%%%%%%
%%%%%%%%%%%%%%%%%%%%%%%%%%%%%%%%%%%%%%%%%%%%%%%%%%%%%%%%%%%%%%%%%%%%%%%%%%%%%%%%%%%%%%%%%%%%%%%%%%%%%
%%%%%%%%%%%%%%%%%%%%%%%%%%%%%%%%%%%%%%%%%%%%%%%%%%%%%%%%%%%%%%%%%%%%%%%%%%%%%%%%%%%%%%%%%%%%%%%%%%%%%

\setcounter{equation}{0}
\section{Introduction}
In this paper we apply compensated convex transforms \cite{Zha08,ZOC14a,ZOC14b,ZCO15} to define
Lipschitz continuous and smooth ($C^{1,1}$)  geometric approximations and interpolations
for bounded real-valued functions sampled from either a compact set $K$ in $\mathbb{R}^n$ or the
complement $K=\mathbb{R}^n\setminus \Omega$ of a bounded open set $\Omega$.
The former is motivated by approximating or interpolating sparse data or contour lines and the latter  by the so-called inpainting problem in image processing \cite{CS05}, where some parts of the image content are missing and the aim is to use other parts of the image to repair or reconstruct the missing parts. We first define two one-sided
approximations, called  upper and lower approximations, from above and below the graph of the sampled function respectively, and
then an average approximation. By using mixed compensated convex transforms \cite{Zha08}, we  will also define a smooth ($C^{1,1}$)
average approximation.
Our central aim here is to develop a mathematical theory for  these average approximations.
 Applications of this theory to  level-set reconstruction, scattered data interpolation and inpainting  will be presented, together with some prototype examples, in a follow-on paper \cite{ZOC15}. 

\medskip
\noindent Before relating our  results to previous work on approximations and interpolations of sampled functions,
we first recall the notions of  quadratic compensated convex transforms
of bounded functions  and present our definitions of upper, lower and average approximations (note that compensated convex transforms can be defined under more general growth conditions than those given here \cite{Zha08}).

\medskip
\noindent Suppose $f:\mathbb{R}^n\mapsto \mathbb{R}$ is bounded.
The quadratic lower and upper compensated convex transform \cite{Zha08}
(lower and upper transforms for short)  are defined for each $\lambda>0$ by
\begin{equation}\label{Eq.Def.UpLwTr}
	C^l_\lambda(f)(x)=\co[\lambda|\cdot|^2+f](x)-\lambda|x|^2,\quad\text{resp.}\quad
	C^u_\lambda(f)(x)=\lambda|x|^2-\co[\lambda|\cdot|^2-f](x),\qquad x\in\mathbb{R}^n,
\end{equation}
where $|x|$ is the standard Euclidean norm of $x\in \mathbb{R}^n$ and $\co[g]$ denotes the convex
envelope \cite{H-UL01,Roc70} of a function $g:\mathbb{R}^n\mapsto \mathbb{R}$ that is bounded below.
For given $\lambda>0$ and $\tau>0$, two quadratic mixed compensated convex transforms \cite{Zha08}
(mixed transforms for short)
are defined, respectively, by $ C^u_\tau(C^l_\lambda(f))$ and $C^l_\tau(C^u_\lambda(f))$.

\medskip
\noindent One key property of the compensated convex transforms, established in \cite{Zha08}, is that $C^l_\lambda(f)$
(respectively, $C^u_\lambda(f)$) realises a `tight' approximation of $f$ from below (respectively, from above),
in the sense that if $f$ is $C^{1,1}$ in a neighbourhood of some $x_0$, then there is  a finite $\Lambda>0$, such
that $f(x_0)=C^l_\lambda(f)(x_0)$ (respectively, $f(x_0)=C^u_\lambda(f)(x_0)$) whenever $\lambda\geq \Lambda$.
A second important property is that of locality.  
Since the definitions\eqref{Eq.Def.UpLwTr} involve the evaluation of the convex envelope
of functions \cite{H-UL01,Roc70}, one might think that these notions are global in nature, that is,
the values of these transforms at a given point might involve values  of the original function far away
from the point. However, the locality property for compensated convex transforms \cite[Theorem 3.10]{ZOC14a} states
that if $f$ is bounded, i.e., $|f(x)|\leq M$ in $\mathbb{R}^n$ for some $M>0$, then
the values of $C^l_\lambda(f)(x_0)$ and $C^u_\lambda(f)(x_0)$ depend only on the values of $f$ in
the closed ball $\bar{B}(x_0;\,R)$ with $R=2\sqrt{2}\sqrt{M/\lambda}$.
As a result, these apparently global transforms are, in fact, local.

\medskip
\noindent In this paper, we mainly consider two types of data sets in $\mathbb{R}^n$, given that
the typical applications  we have in mind  are approximation
of sparse data and of contour lines, and inpainting of damaged images.
We therefore assume in the following that, unless otherwise specified,
$K\subset\mathbb{R}^n$ is
either a compact set or the complement of a bounded open
set $\Omega\subset\mathbb{R}^n$, i.e. $K=\mathbb{R}^n\setminus \Omega$. We denote  by
$f:\R^n\mapsto \R$ the underlying function to be approximated. The function $f_K:K\subset \R^n\mapsto \mathbb{R}$ is
our sampled function defined by $f_K(x) = f(x)$ for $x \in K$,  and $\Gamma_{f_K}:=\{(x,f_K(x)),\, x\in K\}$ is its graph.

\medskip
\noindent  Let $K\subset \R^n$  be a non-empty closed set and suppose that for some constant
$A_0>0$, $|f_K(x)|\leq A_0$ for  all $x\in K$.
Given $M>0$, we define two functions extending $f_K$ to $\R^n\setminus K$, namely
\begin{equation}\label{Eq.Def.ExtFnct}
\begin{array}{lll}
	\displaystyle f^{-M}_K(x) & \displaystyle = f(x)\chi_K(x)-M\chi_{\R^n\setminus K}
				& \displaystyle =
		\left\{\begin{array}{ll}
				f_K(x),	& x\in K,\\[1.5ex]
				-M,	& x\in\mathbb{R}^n\setminus K\,;
			\end{array}\right.\\[2.5ex]
	\displaystyle f^{M}_K(x) & \displaystyle = f(x)\chi_K(x)+M\chi_{\R^n\setminus K}
				& \displaystyle =
		\left\{\begin{array}{ll}
			f_K(x),	& x\in K,\\[1.5ex]
			M,	& x\in\mathbb{R}^n\setminus K\,,
		\end{array}\right.
\end{array}
\end{equation}
\noindent where  $\chi_G$ denotes the characteristic
function of a set $G$.

\medskip
\begin{defi}\label{Def.LwrUpAv}
	For $M>0$, the {\bf upper compensated convex approximation} with
	scale $\lambda>0$ for the sampled function $f_K:K \to \R$ is defined by
\begin{equation}\label{Eq.Def.UpAprx}
	U^{M}_\lambda(f_K)(x)= C^u_\lambda(f^{-M}_K)(x),\quad x\in \mathbb{R}^n\,.
\end{equation}

\medskip
\noindent The {\bf lower compensated convex approximation}  with
scale $\lambda>0$ for the sampled function $f_K:K \to \R$ is defined by
\begin{equation}\label{Eq.Def.LwAprx}
	L^{M}_\lambda(f_K)(x)= C^l_\lambda(f^{M}_K)(x),\quad x\in \mathbb{R}^n\,.
\end{equation}

\medskip
\noindent The {\bf average compensated convex approximation}  with scale $\lambda>0$
for the sampled function $f_K:K \to \R$ is defined by
\begin{equation}\label{Eq.Def.AvAprx}
	A^{M}_{\lambda}(f_K)(x)= \frac{1}{2}
		\left(C^l_\lambda(f^{M}_K)(x)+C^u_\lambda(f^{-M}_K)(x)\right),
		\quad x\in \mathbb{R}^n.
\end{equation}

\medskip
\noindent The {\bf mixed average compensated convex approximation} with scales $\lambda>0$
and $\tau>0$ for the sampled function $f_K:K \to \R$ is defined by
\begin{equation}
	(SA)^{M}_{\tau,\lambda}(f_K)(x)=
	\frac{1}{2}( C^u_\tau(C^l_\lambda(f_K^M))(x)+C^l_\tau(C^u_\lambda(f_K^{-M}))(x)\,,
	\quad x\in \mathbb{R}^n\,.
\end{equation}
\end{defi}

\noindent In the following, we  refer to the approximations in Definition \ref{Def.LwrUpAv}, for short, as the  upper,
lower, average and mixed approximations.

\medskip
\noindent Note that since the
mixed compensated convex transforms are $C^{1,1}$ functions \cite[Theorem 2.1(iv) and Theorem 4.1(ii)]{Zha08},
the mixed average approximation $(SA)^{M}_{\tau,\lambda}$ is a smooth version of our average approximation.
Also, for a bounded function $f:\mathbb{R}^n\mapsto\mathbb{R}$,
satisfying $|f(x)|\leq M$, $x\in\mathbb{R}^n$ for some constant $M>0$,
we have the following estimates \cite[Theorem 3.13]{ZOC14a}
\begin{equation*}
	0\leq C^u_\tau(C^l_\lambda(f))(x)-C^l_\lambda(f)(x)\leq \frac{16M\lambda}{\tau},\quad
	0\leq C^u_\lambda(f)(x)-C^l_\tau(C^u_\lambda(f))(x)\leq \frac{16M\lambda}{\tau}
\end{equation*}
for all $x\in\mathbb{R}^n$, $\lambda>0$ and $\tau>0$, and hence can easily show that for any closed set $K\subset \R^n$, 
\begin{equation*}
	|(SA)^M_{\tau,\lambda}(f_K)(x)-A_\lambda^M(f_K)(x)|\leq \frac{16M\lambda}{\tau},
	\quad  x\in\R^n\,.
\end{equation*}
This  implies that for given $\lambda>0$ and $M>0$,  the mixed approximation $(SA)^M_{\tau,\lambda}(f_K)$
converges to the basic average approximation $A_\lambda^M(f_K)$ uniformly in $\mathbb{R}^n$ as $\tau\to\infty$,
with rate of convergence $16M \lambda/\tau$.

\begin{nota} \label{Rmk.AprxSet}
We can additionally consider the families of average approximations 
\[
	A_{\lambda,s}^M(f_K)(x)=sC^l_\lambda(f^M_K)(x)+(1-s)C^u_\lambda(f^{-M}_K)(x),\quad s\in [0,\,1]
\]
and
\[
	(SA)_{\tau,\lambda,s}^M(f_K)(x)=sC^u_\tau(C^l_\lambda(f^M_K))(x)+
				(1-s)C^l_\tau(C^u_\lambda(f^{-M}_K))(x),\quad s\in [0,\,1]\,.
\]
These more general average approximations give  some flexibility when dealing with sets
which are not  graphs of single-valued functions.
For instance, suppose $X\subset \mathbb{R}^n\times \mathbb{R}$ is a finite set.
Let $K=\mathbb{P}_{\mathbb{R}^n}(X):=\{x_1,\ldots,x_n\}$ be the orthogonal projection of $X$ to $\mathbb{R}^n$,
and  for $x \in K$, define
\[
	\check{f}_K(x)=\inf\{v,  (x,v)\in X\},\qquad \hat{f}_K(x)=
	\sup\{v, \;  (x,v)\in X\}\,.
\]
Then $\check{f}_K(x)\leq \hat{f}_K(x)$ and $\check{f}_K$, $\hat{f}_K$  are both single-valued functions.
We can then define
\begin{equation}
	A_{\lambda,s}^M(X)(x): =sC^l_\lambda(\check{f}^M_K)(x)+(1-s)C^u_\lambda(\hat{f}^{-M}_K)(x)
\end{equation}
for suitable $M$, and optimise $A_{\lambda,s}^M(X)$ with respect to $s\in [0,\,1]$ to find
a good approximation of the set $X$ by a single-valued function. For example, we may consider 
the following nonlinear least square approximation of the data set by the family of
functions $A^M_{\lambda,s}(X)$,
\[
	\inf_{s\in[0,\,1]}\sum_{i=1}^n\max\left\{|A^M_{\lambda,s}(X)(x_i)-v|^2,\;(x_i,v)\in X \right\}\,.
\]
However, we do not explore this further here,  instead focussing on our  basic
average approximation  $A_{\lambda}^M(f_K)$ and the mixed approximation $(SA)^{M}_{\tau,\lambda}$.
\end{nota}

\noindent If we consider the special case where $K$ is a finite set, the average approximation
$A_{\lambda}^M(f_K)$ defines an approximation  for the scattered data
$\Gamma_{f_K}=\{(x,f_K(x)),\; x\in K\}$. Moreover, although our extended
functions are defined in the whole space $\mathbb{R}^n$, when $K$ is compact
we are interested  only in the values of our average approximation
$A_{\lambda}^M(f_K)(x)$ for $x$ in the convex hull $\co[K]$ of the sampled set $K$.
If $K$ is the complement of a bounded open set $\Omega\subset\mathbb{R}^n$,
we will consider the values of $A_{\lambda}^M(f_K)(x)$ for $x$ in the whole space $\mathbb{R}^n$
or in a large domain containing $\bar\Omega$.

\medskip
\noindent Theoretically, we may also set $M=+\infty$ and consider the following
functions, which are commonly used in convex analysis, in place of \eqref{Eq.Def.ExtFnct}:
\begin{equation}
	f^{-\infty}_K(x)=\left\{\begin{array}{l} f(x),\quad x\in K,\\
	-\infty,\quad x\in\mathbb{R}^n\setminus K;\end{array}\right.\qquad
	f^{+\infty}_K(x)=\left\{\begin{array}{l} f(x),\quad x\in K,\\
	+\infty,\quad x\in\mathbb{R}^n\setminus K.\end{array}\right.
\end{equation}
This method of extension can help to establish better approximation results than those obtained using $f^{-M}_K$ and $f^M_K$ 
(compare Theorem \ref{Thm.AprxUF.Cmpct} with Theorem \ref{Thm.AprxCnt}). 
Note, however, that 
 the corresponding average
approximation, $$A_{\lambda}^{\infty}(f_K)(x) : = \frac{1}{2}
		\left(C^l_\lambda(f^{+\infty}_K)(x)+C^u_\lambda(f^{-\infty}_K)(x)\right),
		\quad x\in \mathbb{R}^n,$$ 
is not Hausdorff stable with respect to sample sets in the sense introduced in Section \ref{Sec.HauStb}, in contrast to the basic average approximation $A^M_{\lambda}(f_K)$.

\bigskip
\noindent We turn now to some  background and motivation. Selected recent developments on approximation and interpolation
methods are discussed in \cite{JBHSS06}. The literature on approximation and interpolation theory for sampled
functions over the real line $\mathbb{R}$ by polynomials and
other functions is quite rich  \cite{Tim94,CL00}. When $n\geq 2$, however, many methods for $\mathbb{R}$ no longer apply directly to
$\mathbb{R}^n$. In particular, there is no  direct construction of interpolations for randomly placed
sample functions based on one-dimensional interpolation methods.
For scattered data, Delaunay triangulation-based direct spline designs have been
widely used in computational geometry \cite{OBSC00,Ede87}.
Thin plate spline methods,  variational methods,
which are related to  radial basis function methods, and more
general kernel methods, have been used extensively
in applications \cite{Wen05,Buh04,JBHSS06}, and
morphological reconstruction methods based on level
sets  using geodesic distance \cite[see Ch.6.4 and Ch.7.1.2]{Soi04} have also been
developed.
 Nonlinear partial differential equations
and variational methods using various total variation (TV) based models \cite{ROF92}
have  been used for image reconstruction problems,
 salt \& pepper noise
reduction \cite{CHN05} and  image inpainting \cite{BSCB00,CS05}. 
Although there  is a  well-developed mathematical theory on the existence and uniqueness
of their weak solutions  \cite{GMS79,ABCM01a,ABCM01b,BCN02}, the quantitative effectiveness
of such methods  is mostly assessed on the basis of
numerical experiments. 

\medskip
\noindent Note that many methods are sensitive to the type of data to be interpolated or approximated, that is, to the sample type.
The spline function interpolation and finite element based methods design
interpolations require  precise knowledge of the sample locations.
In this case,  Delaunay triangulation or other types of decomposition typically must  be constructed first 
 \cite{OBSC00}. The radial basis function method relies on solving systems
of linear equations \cite{Wen05}. In order to apply any of these methods to interpolate or  approximate
data sets, one has to assume that the data set is scattered, that is, the set is finite
and the points are isolated. If the data set is given by contour lines
(or by level sets), further discretisation is required before such methods can be used.

\medskip
\noindent Our approach, on the other hand, is not sensitive to data types.
We only assume the underlying  function to be bounded and uniformly continuous, and the sample
sets to be compact or to be the complement of a bounded open set. 
In the digital setting, the data  are always
finite sets, and in a `point cloud', a line can be formed by discrete points next to
each other, which, by definition, should not be thought of as scattered data.
Therefore further down sampling might be needed in order to apply spline or radial basis
function methods. But our average approximation $A^M_\lambda(f_K)$, on the contrary, applies directly to these data sets.
In addition,  collected data are  bounded in a given window, and thus the assumption 
of boundedness of the underlying functions covers most
situations in applications. It should be noted that the idea of using averages for 
approximations is  natural 
and has also been used before by several authors, for example \cite{BGLW08,BLT08,BLW07} introduce the notion 
of proximal average, a parametrized convex function that provides a continuous transformation of 
a convex function into another. In \cite{Har09}
this transformation has also been applied to non-convex functions and with non-quadratic weights by exploiting 
its relationship with the Moreau envelopes, and has been used as a fundamental tool to justify 
the application of parallel proximal algorithms in nonsmooth optimization \cite{Yu13,PB13}.

\medskip
\noindent The exact form of an interpolation is of interest but is often not known.  
An advantage of Delaunay triangulation-based spline interpolation methods is that for
simple geometric examples, one can describe precisely what the interpolation is, in contrast to, for instance, radial basis function and partial differential equation based methods.
Although we do not deliberately design the form of our interpolations, it can be shown that our average approximation  $A^M_\lambda(f_K)$ produces
particular forms for us automatically. For example, if $K$ is finite and $\lambda>0$, $M>0$ are large, we will prove in a follow-on paper \cite{ZOC15} that
$A^M_\lambda(f_K)(x)$ is a piecewise affine interpolation from $K$ to $\co[K]$.
We can also give explicit calculations of our approximations in some other simple geometric cases.

\medskip
\noindent A further natural and practical question in data approximation and interpolation is  the stability of a given method. For approximations and interpolations of sampled
functions, we would like to know, for two sample sets which are `close' to each other,
say, under the Hausdorff distance \cite{AT04}, whether the corresponding approximations are close
to each other. It is easy to see that differentiation and integration based approximation
methods are not Hausdorff stable because continuous functions can be sampled over a finite dense set.
One of the advantages of our method is that for a bounded uniformly continuous function $f$,
and for fixed $M>0$ and $\lambda>0$, the mapping $K\mapsto A^M_\lambda(f_K)(x)$ is continuous
with respect to the Hausdorff distance
for compact sets $K$, and the continuity is uniform with respect to $x\in \mathbb{R}^n$.
This means that if another sampled subset $E\subset \mathbb{R}^n$ (finite or compact) is close
to $K$, then the output  $A^M_\lambda(f_E)(x)$ is close to $A^M_\lambda(f_K)(x)$ uniformly with
respect to $x\in \mathbb{R}^n$.
As far as we know, not many known interpolation/approximation methods  share
such a property.

\medskip
\noindent To shed further light on the theory we develop, it is worth observing the connection between the
compensated convex transforms and our proposed average approximation on the one hand, and the
critical mixed Moreau envelopes and  mathematical morphology on the other hand.
The lower and upper transforms can be viewed as `one-step' morphological
opening and closing, respectively \cite{ZOC14a}.
They in fact coincide with the critical mixed Moreau envelopes, that is,
\begin{equation}\label{Eq.MixMor}
	C^l_\lambda(f)(x)=M^\lambda(M_\lambda(f))(x)
	\quad\text{and}\quad C^u_\lambda(f)(x)=M_\lambda(M^\lambda(f))(x)\,,
\end{equation}
where
\begin{equation*}\label{Eq.Def.MixMor}
	M_{\lambda}(f)(x)=\inf\{f(y)+\lambda|x-y|^2,\; y\in\mathbb{R}^n\}
	\quad\text{and}\quad
	M^{\lambda}(f)(x)=\sup\{f(y)-\lambda|x-y|^2,\; y\in\mathbb{R}^n\}
\end{equation*}
are the lower and upper Moreau envelopes \cite{Mor65,Mor66,LL86,AA93,CS04}, respectively. If we denote by
$b_{\lambda}(x)=-\lambda|x|^2$ the quadratic structuring function, introduced for the first time
in \cite{Jac92,BH00,BS92,JD96}, then with the notation
of \cite{Ser82,Soi04}, we have\footnote{In convex analysis, the infimal convolution of $f$ with $g$
is denoted in  \cite{Roc70} 
as $f\Box g$ and is defined as $(f\Box g)(x)=\underset{y}{\inf}\, \{f(y)+g(x-y)\}$, 
whereas in \cite{Hir94} the deconvolution of $f$ with $g$ is denoted as $f\ominus g$ and, 
under the condition that for some $x_0\in \mathbb{R}^n$ and $r\in\mathbb{R}$, we have $f(x)\leq g(x-x_0)+r$ 
for all $x\in\mathbb{R}^n$, is defined as $(f\ominus g)(x)=\underset{y}{\sup}\, \{f(x-y)-g(y)\}$. Thus 
$M_{\lambda}(f)$ is the inf-convolution of $f$ with $-b_{\lambda}$ whereas
$M^{\lambda}(f)$ is the deconvolution of $f$ with $b_{\lambda}$.
}
\begin{equation*}
\begin{split}
	M_{\lambda}(f)(x)&=\inf\{f(y)-b_{\lambda}(x-y),\; y\in\mathbb{R}^n\}=f\ominus b_{\lambda}\,,\\[1.5ex]
	M^{\lambda}(f)(x)&=\sup\{f(y)+b_{\lambda}(x-y),\; y\in\mathbb{R}^n\}=f\oplus  b_{\lambda}
\end{split}
\end{equation*}
that is, the Moreau lower and upper envelopes can be viewed as `greyscale' erosion and dilation with quadratic
structuring function, respectively \cite{BH00,Luc09}. Compared with \eqref{Eq.MixMor}, we thus have 
\begin{equation*}
	C_{\lambda}^l(f)=(f\ominus b_{\lambda})\oplus b_{\lambda}
	\quad\text{and}\quad
	C_{\lambda}^u(f)=(f\oplus b_{\lambda}) \ominus b_{\lambda}\,,
\end{equation*}
and hence, using the definition  of $A^M_{\lambda}(f_K)$, it follows that
\begin{equation*}
	A^M_{\lambda}(f_K)=\frac{1}{2}\left((f^{-M}_K\oplus b_{\lambda})\ominus b_{\lambda} +
	(f^{M}_K \ominus b_{\lambda})\oplus b_{\lambda}\right).
\end{equation*}
Given such an interpretation for $A^M_{\lambda}(f_K)$,
the properties of $A^M_{\lambda}(f_K)$ could therefore also be analysed by tools
from the theory of morphological filtering \cite{SV92,GMM04}.

\medskip
\noindent The plan of the rest of this paper is as follows. In  Section \ref{Sec.Nots},
we introduce notation and recall some useful results from convex analysis
and compensated convexity theory.
Our main sources of references for convex analysis are \cite{H-UL01,Roc70,RW98} whereas for the
properties of compensated convex transforms we refer to \cite{Zha08,ZOC14a,ZOC14b,ZCO15}.
In Section \ref{Sec.IntAprx} we state our general approximation/interpolation
theorems for a compact sample set $K\subset \mathbb{R}^n$ (Theorem \ref{Thm.AprxUF.Cmpct}) and for
$K=\mathbb{R}^n\setminus \Omega$ with $\Omega$ a bounded open set (Theorem \ref{Thm.AprxCnt}).
We consider uniformly continuous, Lipschitz and $C^{1,1}$ functions
$f:\mathbb{R}^n\mapsto\mathbb{R}$ as the underlying functions to be approximated.
We show that when $M>0$ is sufficiently large,
$A^M_\lambda(f_K)$ approaches $f_K$ in $K$ as $\lambda \to \infty$.
If $f$ is a $C^{1,1}$ function, we also show that  $A^M_\lambda(f_K)$ is an
interpolation of $f$ in the convex hull $\co[K]$ of $K$ when $\lambda>0$ is large enough.
For points $x$ in $\co[K]\setminus K$, we introduce the notion of convex density radius $r_c(x)$
which is the smallest radius of a closed ball $\bar{B}(x;\,r_c(x))$ such that $x$
is in the convex hull of $K\cap \bar{B}(x;\,r_c(x))$.
We  use  $r_c(x)$ to bound the errors of our approximations
$|A^\infty_\lambda(f_K)(x)-f(x)|$.
For a finite $M>0$ and for a compact sample set $K$, we extend $f_K$ to be a constant $c_0$ outside
a large ball $B(0;\,R)$ containing $K$ and define $K_R=K\cup B^c(0;\,R)$. We then prove
similar error estimates to those obtained for $A^\infty_\lambda(f_K)-f$ (Corollary \ref{Cor.AprxCnt}). 
For example, for a bounded uniformly continuous function $f$,
\begin{equation*}
	|A^M_\lambda(f_{K_R})(x)-f(x)|\leq \omega(r_c(x)+a/\lambda+\sqrt{2b/\lambda}),\quad x\in\mathbb{R}^n\,,
\end{equation*}
where $\omega:[0,\,+\infty)\mapsto [0,\,+\infty)$ is the least concave majorant of the
modulus of continuity of
the function $f$ \cite{{DL93}} which satisfies $\omega(t)\leq at+b$ for $t\geq 0$ and some
constants $a>0$ and $b\geq 0$.
Better estimates  are also established for Lipschitz functions and for $C^{1,1}$ functions.

\medskip
\noindent In Section \ref{Sec.HauStb} we state the Hausdorff stability property for the
average approximation $A^M_\lambda(f_K)$ of bounded uniformly continuous
functions, showing that
given two non-empty  closed sets $K$ and $E$, for fixed $M>0$ and $\lambda>0$,
 $|A^M_\lambda(f_{K})(x)-A^M_\lambda(f_{E})(x)|$ is uniformly small
in $\mathbb{R}^n$ with explicit estimates whenever $K$ and $E$ are closed.
For a bounded Lipschitz function
$f:\mathbb{R}^n\mapsto\mathbb{R}$ with $|f(x)|\leq A_0<M$ for some constant $A_0>0$
and  all $x\in \R^n$,
the mapping $K\mapsto A^M_\lambda(f_{K})(x)$
is Lipschitz continuous with respect to the Hausdorff metric, uniformly with respect to
$x\in \mathbb{R}^n$. This result generalises an earlier Hausdorff-Lipschitz continuity result
for the upper transform of characteristic functions
$K\mapsto C^u_\lambda(\chi_K)$ established in \cite[Theorem 5.5]{ZOC14a}.
We conclude Section \ref{Sec.HauStb} by proving  regularity properties of our approximations.
For example, we show that $A^M_\lambda(f_{K})$ is a globally Lipschitz function
in $\mathbb{R}^n$ and give an explicit estimate of its Lipschitz constant.

\medskip
\noindent The proofs of  our main results are presented in Section \ref{Sec.Prfs}.

\medskip
\noindent In the follow-on paper \cite{ZOC15} we will present some applications of the theory developed here, such
as interpolation and approximation of scattered data and for contour lines.
We will also give some prototype examples with analytical expressions of our approximations,
and numerical experiments on salt-and-pepper denoising, inpainting and contour-line based reconstructions.

%%%%%%%%%%%%%%%%%%%%%%%%%%%%%%%%%%%%%%%%%%%%%%%%%%%%%%%%%%%%%%%%%%%%%%%%%%%%%%%%%%%%%%%%%%%%%%%%%%%%
%%%%%%%%%%%%%%%%%%%%%%%%%%%%%%%%%%%%%%%%%%%%%%%%%%%%%%%%%%%%%%%%%%%%%%%%%%%%%%%%%%%%%%%%%%%%%%%%%%%%
%%%%%%%%%%%%%%%%%%%%%%%%%%%%%%%%%%%%%%%%%%%%%%%%%%%%%%%%%%%%%%%%%%%%%%%%%%%%%%%%%%%%%%%%%%%%%%%%%%%%

\setcounter{equation}{0}
\section{Notation and Preliminaries}\label{Sec.Nots}

In this section we collect  basic results and definitions from convex analysis,  referring to 
\cite{H-UL01,Roc70,RW98} for further
references and proofs, and recall the notion of the Hausdorff distance between two non-empty sets.
We then list some selected basic properties of compensated convex transforms
\cite{Zha08,ZOC14a,ZOC14b,ZCO15} that will be needed in the following.

\begin{prop}\label{SecNot.Prp.GlCnvx}
	Let $f:\mathbb{R}^n\mapsto \mathbb{R}$ be coercive in the sense that
	$f(x)/|x|\to\infty$ as $|x|\to \infty$, and $x_0\in\R^n$. Then
   \begin{itemize}
	\item[(i)]
		The value $\co\left[f\right](x_0)$ of the convex envelope of $f$ at
		$x_0\in \mathbb{R}^n$ is given by
		\begin{equation}\label{Eq.Prp.GlCnvx.1}
			 \co\left[f\right](x_0)=\underset{i=1,\ldots,n+1}{\inf}\,\Bigg\{ \sum^{n+1}_{i=1}\lambda_if(x_i):\;
			        \sum^{n+1}_{i=1}\lambda_i=1,\; \sum^{n+1}_{i=1}\lambda_ix_i=x_0,\;
				\lambda_i\geq 0,\; x_i\in \R^n
			\Bigg\}\,.
		\end{equation}
		If, in addition, $f$ is lower semicontinuous, the infimum is  attained by some
		$(\lambda_i^\ast,\, x_i^\ast)$ for $i=1,2,\dots,n+1$ with
		$\left(x_i^\ast, f(x_i^\ast)\right)$'s lying in the intersection
		of a supporting plane of the epigraph of $f$, $\epi(f)$, and
		$\epi(f)$ (see \cite[Lemma 3.3(ii)]{BH96}, \cite[Theorem 2.1]{GR90}, and
		\cite[Corollary 3.47]{RW98}). In this case, 
		\begin{equation}\label{Eq.Prp.GlCnvx.1a}
			\co\left[f\right](x_0)=\sum^{n+1}_{i=1}\,\lambda_i^\ast f(x_i^\ast)\,.
		\end{equation}
	\item[(ii)]
		The value $\co\left[f\right](x_0)$, for $f$ taking only finite values, can also be obtained
		as follows:
		\begin{equation}\label{Eq.Prp.GlCnvx.2}
			\co\left[f\right](x_0)=\sup\left\{
				\ell(x_0):\; \ell\;\;\text{\rm affine}\quad\text{\rm and}\quad
				\ell(y)\leq f(y)\;\;\text{\rm for all }y\in\R^n
				\right \}
		\end{equation}
	with the $\sup$ attained by an affine function $\ell^\ast\in\Aff(\R^n)$.
	\item[(iii)] If $f$ is differentiable at $x_0$ and
		\begin{equation}\label{Eq.Prp.GlCnvx.3}
			f(x)\geq f(x_0)+Df(x_0)\cdot(x-x_0)\quad \text{for all }x\in \R^n\,,
		\end{equation}
		then $\co[f](x_0)=f(x_0)$.
     \end{itemize}
\end{prop}

\medskip
Next we recall the definition of Hausdorff distance between two
non-empty sets \cite{AT04}, which measures how far the  sets are from each other. To do so, we first
need the notion of $\delta$-neighbourhood of a set, and also define the diameter of a set.

\begin{defi}\label{SecNot.Def.NeighSet}
	Given a non-empty subset $E\subset \mathbb{R}^n$ and $\delta>0$, we define the
	$\delta$-neighbourhood $E^\delta$ of $E$ by
	\[
		E^\delta=\{x\in \mathbb{R}^n,\; \dist(x;\, E)<\delta\}\,.
	\]
where $\dist(x;\, E)=\inf\{|x-y|,\,y\in E\}$, and the diameter of $E$  by
$$d_E:=\sup\{|x-y|,\; x,\;y\in E\}.$$

\end{defi}

\medskip

\begin{defi} \label{Sec2.Def.HausDist}
	Let $E,\, F$ be non-empty subsets of $\R^n$. The Hausdorff distance
	between $E$ and $F$ is defined by
	\begin{equation}\label{Sec2.Eq.HausDist}
		\dist_{\mathcal{H}}(E,F)=\inf\left\{\delta>0:
		F\subset E^\delta\; \text{and}\;\, E\subset F^\delta\right\}.
	\end{equation}
\end{defi}

\noindent For general closed sets $K,\, G\subset \mathbb{R}^n$, if there is some $\delta>0$ such that
$K\subset G^\delta,\; G\subset K^\delta$, then the Hausdorff distance between $F$ and $G$ is finite and
is given by \eqref{Sec2.Eq.HausDist}. Otherwise we say that $\dist_{\mathcal{H}}(K,G)=+\infty$.

\bigskip
\noindent We now list some properties of the quadratic compensated convex transforms.
Recall first the following ordering properties  \cite{Zha08}:

\begin{equation}\label{Sec2.Eq.Ordering}
	C^l_\lambda(f)(x)\leq f(x)\leq C^u_\lambda(f)(x),\quad x\in\mathbb{R}^n\,,
\end{equation}

\noindent whereas for $f\leq g$ in $\mathbb{R}^n$, we have that

\begin{equation}\label{Sec2.Eq.Ineq}
	C^l_\lambda(f)(x)\leq C^l_\lambda(g)(x)\quad\text{and}\quad
	C^u_\lambda(f)(x)\leq C^u_\lambda(g)(x),\quad x\in\mathbb{R}^n\,.
\end{equation}

\noindent Furthermore, the compensated convex transforms are affine invariant \cite{ZOC14a}, that is,
\begin{equation}
	C_{\lambda}^{l}(f+\ell)=C_{\lambda}^{l}(f)+\ell
	\quad\text{and}\quad
	C_{\lambda}^{u}(f+\ell)=C_{\lambda}^{u}(f)+\ell
\end{equation}
\noindent where $\ell$ is any affine function, and we also have \cite[Theorem 2.1(iii)]{Zha08}
\begin{equation}
	C_{\tau}^{u}(C_{\lambda}^{u}(f))=\left\{\begin{array}{ll}
			\displaystyle C_{\lambda}^{u}(f) & \displaystyle  \text{if }\tau\geq\lambda\,,\\[1.5ex]
			\displaystyle C_{\tau}^{u}(f)	 & \displaystyle  \text{if }\tau\leq\lambda\,;
					\end{array}\right.
					\quad\text{and}\quad
	C_{\tau}^{l}(C_{\lambda}^{l}(f))=\left\{\begin{array}{ll}
			\displaystyle C_{\lambda}^{l}(f) & \displaystyle  \text{if }\tau\geq\lambda\,,\\[1.5ex]
			\displaystyle C_{\tau}^{l}(f)	 & \displaystyle  \text{if }\tau\leq\lambda\,.
					\end{array}\right.
\end{equation}

\bigskip
\noindent The following translation-invariance property will  often be used in our proofs, since it allows us to refer our arguments to the point
$x_0=0$ without loss of generality.

\begin{prop}\label{Prp.Tra} (Translation-invariance property)
	For any $f:\mathbb{R}^n\mapsto \mathbb{R}$
	bounded below and for any affine function $\ell:\mathbb{R}^n\mapsto \mathbb{R}$,
	$\co[f+\ell]=\co[f]+\ell$. Consequently,
	both $C^u_\lambda(f)$ and $C^l_\lambda(f)$ are translation invariant against
	the weight function, that is:
	\begin{equation*}\label{Eq.Prp.Tra.1}
	\begin{split}
		&C^l_\lambda(f)(x)=\co\left[\lambda |(\cdot)-x_0|^2+f\right](x)-\lambda |x-x_0|^2\,,\\[1.5ex]
		&C^u_\lambda(f)(x)=\lambda |x-x_0|^2-\co\left[\lambda |(\cdot)-x_0|^2-f\right](x)
	\end{split}
	\end{equation*}
	for all $x\in \R^n$ and for every fixed $x_0$. Hence, at $x_0$,
	\begin{equation*}\label{Eq.Prp.Tra.2}
			C^l_\lambda(f)(x_0)=\co[\lambda  |(\cdot)-x_0|^2+f](x_0)\,, \quad
			C^u_\lambda(f)(x_0)=-\co[\lambda |(\cdot)-x_0|^2-f](x_0)\,.
	\end{equation*}
\end{prop}

\noindent For some theoretical developments and proofs, it can  be convenient to view the
lower and upper compensated convex transforms as parametrized semiconvex and
semiconcave envelopes, respectively. We recall the following definition from   \cite{CS04,CSW95}.

\begin{defi}\label{Def.Semiconcave}
	A function $f:\mathbb{R}^n\mapsto\mathbb{R}$ is called $2\lambda$-semiconvex
	(respectively, $2\lambda$-semiconcave) if $x\mapsto f(x)+\lambda|x|^2$
	(respectively, if $x\mapsto -f(x)+\lambda|x|^2$) is convex.
\end{defi}	

\begin{nota}
In convex analysis, the semiconvexity property as given by Definition \ref{Def.Semiconcave} 
is sometimes also referred to as the
uniform lower-$C^2$ property; compare Definition \ref{Def.Semiconcave} with that of lower-$C^2$ in 
\cite[page 228]{Bou12}. Such functions enjoy local regularity properties; note, for instance,
the characterization of the locally Lipschitz functions as locally lower-$C^2$ - see \cite[Theorem 6]{Roc82}
and \cite[Theorem 2.1.7]{CS04}.
\end{nota}

\medskip
\noindent In our approximation theorems for bounded and uniformly continuous functions $f$,
we make use of the modulus of continuity of $f$, which
 is defined as follows \cite{DL93}.
\begin{defi}\label{Sec2.Def.MoC}
	Let $f:\R^n \mapsto \R$ be a bounded and uniformly continuous function in $\R^n$. Then,
	\begin{equation}\label{Sec2.Eq.Def.MoC}
		\omega_f:t\in[0,\,\infty)\mapsto
		\omega_{f}(t)=\sup\Big\{|f(x)-f(y)|:\,x,y\in\R^n\text{ and }|x-y|\leq t\Big\}
	\end{equation}
	is called the modulus of continuity of $f$.
\end{defi}

\noindent The modulus of continuity of $f$ has the following properties.

\begin{prop}\label{Sec2.Pro.MoC}
Let $f:\R^n \mapsto \R$ be a bounded and uniformly continuous function in $\R^n$.
Then the modulus of continuity $\omega_f$ of $f$ satisfies the following properties:
\begin{equation}\label{Sec2.Eq.Pro.MoC}
	\begin{array}{ll}
		$(i)$	&\omega_f(t)\to \omega_f(0)=0,\text{ as }t\to 0;	\\[1.5ex]
		$(ii)$	&\omega_f \text{ is non-negative and non-decreasing continuous function on }[0,\infty);	\\[1.5ex]
		$(iii)$	&\omega_f \text{ is subadditive: }\omega_f(t_1+t_2)\leq\omega_f(t_1)+\omega_f(t_2)
			\text{ for all }t_1,\,t_2\geq 0\,.
	\end{array}
\end{equation}
\end{prop}

\medskip

\noindent Any function $\omega$ defined on $[0,\,\infty)$ and satisfying \eqref{Sec2.Eq.Pro.MoC}{\it (i), (ii), (iii)} is 
called {\it a
modulus of continuity}. A modulus of continuity
$\omega$ can be bounded from above by an affine function (see \cite[Lemma 6.1]{DL93}),
that is, there exist  constants $a>0$ and $b\geq 0$ such that
\begin{equation}\label{Sec2.Eq.Pro.MoC.0}
	\omega(t)\leq at+b\quad(\text{for all }t\geq 0).
\end{equation}
As a result, given $\omega_f$,  one can define the least concave majorant of $\omega_f$,
which we denote by $\omega$, which is also
a modulus of continuity with the property (see \cite{DL93})
\begin{equation}\label{Sec2.Eq.Pro.MoC.01}
	\frac{1}{2}\omega(t)\leq \omega_f(t) \leq \omega(t)\quad(\text{for all }t\in[0,\,\infty))\,.
\end{equation}

%%%%%%%%%%%%%%%%%%%%%%%%%%%%%%%%%%%%%%%%%%%%%%%%%%%%%%%%%%%%%%%%%%%%%%%%%%%%%%%%%%%%%%%%%%%%%%%%%%%%
%%%%%%%%%%%%%%%%%%%%%%%%%%%%%%%%%%%%%%%%%%%%%%%%%%%%%%%%%%%%%%%%%%%%%%%%%%%%%%%%%%%%%%%%%%%%%%%%%%%%
%%%%%%%%%%%%%%%%%%%%%%%%%%%%%%%%%%%%%%%%%%%%%%%%%%%%%%%%%%%%%%%%%%%%%%%%%%%%%%%%%%%%%%%%%%%%%%%%%%%%

\setcounter{equation}{0}
\section{Approximations and Interpolations}\label{Sec.IntAprx}

In this section we consider the general theory of our interpolation and approximation
problem when either $K\subset\mathbb{R}^n$ is compact or
$K=\mathbb{R}^n\setminus \Omega$ where $\Omega\subset \mathbb{R}^n$ is a bounded open set.

\medskip
\noindent Let $f:\mathbb{R}^n\mapsto\mathbb{R}$  be a bounded function and denote by 
$f_K:K\subset \R^n\mapsto \mathbb{R}$
the restriction of $f$ to $K$. A function $g:\co[K]\subset \R^n\mapsto \mathbb{R}$ is said to be an 
{\it interpolation} of  $f_K$ if $g=f$ in $K$, while for $\lambda>0$, a family of
functions $g_\lambda:\co[K]\subset \R^n\mapsto \mathbb{R}$ is said to {\it approximate} $f$ 
if $\displaystyle \lim_{\lambda\to+\infty} g_\lambda=f$
uniformly in $K$.

\medskip
\noindent We will see that the precise approximation and interpolation properties of $f_K$  depend on the smoothness
of the function $f$ under consideration. 

\medskip
\noindent The following is a first simple observation.

\begin{prop}\label{Prp.IntSemi}
	Let $f:\mathbb{R}^n\mapsto\mathbb{R}$  be a bounded $2\lambda$-semiconvex
	(respectively, $2\lambda$-semiconcave) function and $K\subset \mathbb{R}^n$
	a non-empty closed set.
	If $|f(x)| < M$ for all $x \in \R^n$,    then for any
	$\tau\geq \lambda$, $C^l_\tau(f^M_K)$ (respectively, $C^u_\tau(f^{-M}_K)$)
	is an interpolation of $f_K$, that is,
	\[
		C^l_\tau(f^M_K)(x)=f(x) \;\;  \mbox{(respectively, $C^u_\tau(f^{-M}_K)(x)=f(x)$),\; } \;\; x\in K.
	\]
\end{prop}

\medskip
\noindent In general, if we are given sample values only in a closed set without any knowledge of the
underlying function, we do not  know whether or not our transforms are approximations of the
original function. However, for any bounded function $f:\mathbb{R}^n\mapsto \mathbb{R}$,
we have \cite[Prop. 3.1]{ZOC14a}
\begin{equation}
\label{Sec3.semicont}
	C^l_\lambda(f)(x)=C^l_\lambda(\underline{f})(x),\quad\text{and}\quad
	C^u_\lambda(f)(x)=C^u_\lambda(\overline{f})(x)
\end{equation}
for all $x\in \mathbb{R}^n$, where  $\overline{f}$ and $\underline{f}$ are the upper and
lower semicontinuous closures of $f$, respectively, and
\begin{equation}
\label{Sec3.convtransforms}
	\lim_{\lambda\to\infty}C^u_\lambda(f)(x) = \overline{f}(x),\qquad
	\lim_{\lambda\to\infty}C^l_\lambda(f)(x) = \underline{f}(x)
\end{equation}
for all $x\in \mathbb{R}^n$. As a result, we have the following general approximation theorem.

\begin{teo}\label{Thm.AprxLimBnd}
Let $f:\mathbb{R}^n\mapsto \mathbb{R}$ be bounded, such that  $|f(x)|<M$ for all $x\in \R^n$, and let $K\subset \mathbb{R}^n$ be non-empty compact set.
Then for all $x\in \mathbb{R}^n$,
\begin{equation} \label{Sec3.Eq.AprxLimBnd}
	\begin{array}{c}
		\displaystyle	\lim_{\lambda\to+\infty} C^u_\lambda(f^{-M}_{K})(x)=\overline{f^{-M}_{K}}(x),\qquad
			\lim_{\lambda\to+\infty}C^l_\lambda(f^M_{K})(x)=\underline{f^M_{K}}(x)\,,
			\\[1.5ex]
		\displaystyle	\lim_{\lambda\to+\infty}A^M_\lambda(f_K)(x)= \frac{1}{2}(\underline{f^{M}_{K}}(x)+
			\overline{f^{-M}_K}(x)),
	\end{array}
\end{equation}
and if $f:\mathbb{R}^n\mapsto\mathbb{R}$ is   continuous,
then for all $x \in K$, 
\begin{equation} \label{Sec3.Eq.AprxLimUC}
	\begin{array}{c}
		\displaystyle \lim_{\lambda\to+\infty} C^u_\lambda(f^{-M}_{K})(x)=f(x),\qquad
				\displaystyle \lim_{\lambda\to+\infty}C^l_\lambda(f^M_{K})(x)=f(x),\\[1.5ex]
		\displaystyle \lim_{\lambda\to+\infty}A^M_\lambda(f_K)(x)= f(x),
	\end{array}
\end{equation}
and the convergence in \eqref{Sec3.Eq.AprxLimUC} is uniform on $K$.
\end{teo}

\noindent Note that the equalities $\overline{f_K^{-M}}= \left(\overline{f}\right)^{-M}_{K}$ and
$\underline{f_{K}^M}=\left(\underline{f}\right)^M_{K}$ do not hold in general.
For example, in  $\mathbb{R}$, if we define $f(x)=1$ if $x$ is rational,
$f(x)=-1$ if $x$ is irrational and take $M=2$, then we have $\overline{f}\equiv 1$ and
$\underline{f}\equiv -1$. But if we choose $K$ to be a finite set of rational numbers,
then $\underline{f_{K}^{M}}=\chi_K+2\chi_{\mathbb{R}\setminus K}$, whereas 
$\left(\underline{f}\right)^M_{K} = -\chi_K+2\chi_{\mathbb{R}\setminus K}$.

\bigskip
\noindent Note also that Theorem \ref{Thm.AprxLimBnd} suggests that we can apply our approximation methods to data sets
which may not define a function,  as discussed in Remark \ref{Rmk.AprxSet}.

\medskip
\noindent The following proposition provides conditions sufficient to
 ensure that our average approximation does not attain the value $M$ or $-M$.

\begin{prop}\label{Prp.StrBndAprx}
Let $K\subset \mathbb{R}^n$ be a non-empty compact set and denote by $d_K$ the diameter of $K$.
Suppose  $f_K:K\subset \R^n\mapsto\mathbb{R}$
is bounded, with $|f_K(x)|\leq A_0$ 
for all $x\in K$. Then for $\lambda>0$ and $M>A_0+\lambda d^2_K$,
\begin{equation*}
	-A_0\leq C^l_{\lambda}(f_K^M)(x)<M,\quad -M<C^u_\lambda(f_K^{-M})(x)\leq A_0,
	\quad -M<A^M_{\lambda}(f_K)(x)<M \,,
\end{equation*}
for all $x\in \co[K]$.
\end{prop}

\noindent Next we state our weak maximum principle. To make our statement simpler,
we  assume that the median of values of $f_K$ over $K$ is zero, which in practice
can be easily satisfied by a simple translation of values.

\begin{teo}\label{Thm.WMP}
Let $K\subset \mathbb{R}^n$ be a non-empty compact set.
Suppose  $f_K:K\mapsto\mathbb{R}$ is bounded and assume that
\begin{equation}
\label{median}
	m(f_K):=\frac{1}{2}\Big(\sup\{f_K(x),\;x\in K\}+\inf\{f_K(x),\;x\in K\}\Big)=0\,.
\end{equation}
Then 
\begin{equation*}
	\inf\{f_K(y),\;y\in K\}\, \leq \, A^\infty_\lambda(f_K)(x) \, \leq \, \sup\{f_K(y),\;y\in K\}
\end{equation*}
for all $x\in \co[K]$ and $\lambda>0$.
\end{teo}

\bigskip
\noindent Before stating the error estimates for our approximations, we
introduce the notions of density radius and convex density radius.

\begin{defi}\label{Sec3.Def.DnstRadius}		
Suppose $K\subset \mathbb{R}^n$ is a non-empty closed
set, and denote by $\dist(x;\,K)$ the Euclidean distance of $x$ to $K$.
\begin{itemize}
	\item[(i)] For $x\in \co[K]$, the density radius $r_d(x)$ of $x$ with respect to $K$ is
		just the Euclidean distance of $x$ to $K$, i.e. we set $r_d(x)=\dist(x;\,K)$,
		whereas the  density radius of $K$ in $\co[K]$
		is defined as
		\[
			r_d(K)=\sup\{r_d(x),\; x\in \co[K]\}\,.
		\]
	\item[(ii)]  For $x\in \co[K]$, consider the balls $B(x;\,r)$ such that
		$x\in \co[\bar{B}(x;\,r)\cap K]$. The convex density radius of $x$
		with respect to $K$ is defined as follows
		\[
			r_c(x)=\inf\{r\geq 0, \; x\in \co[\bar{B}(x;\,r)\cap K]\}\,,
		\]
		whereas the convex density radius of $K$ in $\co[K]$ is defined by
		\[
			r_c(K)=\sup\{r_c(x),\; x\in \co[K]\}\,.
		\]
\end{itemize}
\end{defi}

\noindent From the definition above, we see that if $K$ is compact or $K^c$ is a bounded 
open set,
$r_c(K)$ is finite. 
The convex density radius is zero if $K$ is convex. If $A$ and $B$ are two compact
sets such that $K\subset A\subset B\subset \co[K]$, then
$r_c(A)\geq r_c(B)$. Also, the smaller $r_c(K)$ is, the denser
the set $K$ is in $\co[K]$.
In general, if $K$ is compact, $r_c(K)$ can be as large as the diameter of $K$.
In this case, in order to make the convex density radius of $K$ small,
we require not only the density radius $r_d(K)$ of $K$ in $\co[K]$ to be small
but also that $K$ is `dense' in $\partial \co[K]$, the relative boundary of $\co[K]$.
If $K^c$ is bounded, then $r_c(K)$ can be as large as the diameter of $K^c$.

\medskip
\noindent The following is a simple illustrative example for the case $K$  compact.
Consider the box $D=\{(x,y)\in \mathbb{R}^2,\; |x|< 1,\; |y|<1\}$. For any $\delta>0$,
let $K_0\subset D$ be a finite set with
$r_d(K_0)<\delta$, so that $\bar D\subset K_0^\delta$, and
let $K=K_0\cup\{\pm 1,\pm1\}$.
Then $r_d(K_0)<\delta$, whereas $r_c(K)=1$ if we consider, say,
the point $(1,0)\in \co[K]=\bar{D}$.

\bigskip
\noindent We now  formulate  error estimates for our average approximations.
Consider  first the case when $K$ is compact and $M=+\infty$.
The estimates are expressed in terms of
the modulus of continuity of the underlying uniformly continuous function $f$
and  the convex density radius.
As special cases, we also consider bounded Lipschitz functions and $C^{1,1}$ functions.

\begin{teo}\label{Thm.AprxUF.Cmpct}	
	Suppose $f:\mathbb{R}^n\mapsto \mathbb{R}$ is a bounded uniformly continuous
	function satisfying $|f(x)|\leq A_0$ for some constant $A_0>0$ and all
	$x\in \R^n$, and let $K\subset \mathbb{R}^n$ be a non-empty compact set.
	\begin{itemize}
		\item[(i)] Denote by $\omega$ the least concave majorant of the modulus of continuity $\omega_f$
			   of $f$. Let $a\geq 0$, $b\geq 0$ be such that $\omega(t)\leq at+b$ for $t\geq 0$. 
		      	   Then for all $\lambda >0$ and  $x\in \co[K]$,
			\begin{equation}
				|A^\infty_\lambda(f_K)(x)-f(x)|\leq
				\omega\left(r_c(x)+\frac{a}{\lambda}+\sqrt{\frac{2b}{\lambda}}\right)\,,
			\end{equation}
			where $r_c(x)\geq 0$ is the convex density radius of $x$ with respect to $K$.
		\item[(ii)] If we further assume that  $f$ is a globally Lipschitz function with
				Lipschitz constant $L>0$, then for all $\lambda>0$ and   $x\in \co[K]$,
			\begin{equation}
				|A^\infty_\lambda(f_K)(x)-f(x)|\leq Lr_c(x)+\frac{L^2}{\lambda}\,.
			\end{equation}
		\item[(iii)] If we further assume that  $f$ is a $C^{1,1}$ function satisfying $|Df(x)-Df(y)|\leq L|x-y|$ 
			     for $x,\, y\in\mathbb{R}^n$ and for some fixed $L>0$, then for all $\lambda>L$ and  $x\in \co[K]$,
			     \begin{equation}
					|A^\infty_\lambda(f_K)(x)-f(x)|\leq
					\frac{L}{4}\left(\frac{\lambda+L/2}{\lambda-L/2}+1\right)r_c^2(x).
				\end{equation}
			Furthermore,  in case $(iii)$, $A^\infty_\lambda(f_K)$ is an interpolation of $f^\infty_K$ in $\co[K]$.
	\end{itemize}
\end{teo}

\medskip

\noindent Next we consider the case when $\Omega\subset \mathbb{R}^n$ is a non-empty bounded open set and define
$K=\Omega^c:=\R^n\setminus\Omega$. Clearly,  $\co[K]=\R^n$ for such $K$.
We then have the following estimate of the average approximation $A^M_\lambda(f_{K})$.

\begin{teo}\label{Thm.AprxCnt}
Suppose $f:\mathbb{R}^n\mapsto \mathbb{R}$ is bounded and uniformly continuous, satisfying $|f(x)|\leq A_0$ for
some constant $A_0>0$ and all $x\in \R^n$. Let $\Omega\subset \mathbb{R}^n$ be a bounded open set and $K=\Omega^c$. Denote by
		$d_\Omega$ the diameter of $\Omega$.
\begin{itemize}
	\item[(i)] Let $\omega$ be the least concave majorant of the modulus of continuity $\omega_f$ of $f$.
		Assume $a\geq 0$, $b\geq 0$ are such that $\omega(t)\leq at+b$ for $t\geq 0$.
		Then for $\lambda>0$,  $M>A_0+\lambda d_\Omega^2$
		 and  all $x\in \mathbb{R}^n$, we have
		\begin{equation}
			|A^M_\lambda(f_{K})(x)-f(x)|\leq \omega\left(r_c(x)+\frac{a}{\lambda}+
			\sqrt{\frac{2b}{\lambda}}\right)\,,
		\end{equation}
		where $r_c(x)\geq 0$ is the convex density radius of $x$ with respect to $K$.
	\item[(ii)] If we further assume that $f$ is a globally Lipschitz function with Lipschitz constant $L>0$,
	 then for $\lambda>0$,  $M>A_0+\lambda d_\Omega^2$ and all $x\in \mathbb{R}^n$, we have
	\begin{equation}
		|A^M_\lambda(f_{K})(x)-f(x)|\leq Lr_c(x)+\frac{L^2}{\lambda}\,.
	\end{equation}
	\item[(iii)] If we further assume that $f$ is a $C^{1,1}$ function such that $|Df(x)-Df(y)|\leq L|x-y|$
		for all $x,\,y\in \mathbb{R}^n$, where $L>0$ is a constant,  then for $\lambda>L$, 
		$M>A_0+\lambda d_\Omega^2$ and all $x\in \mathbb{R}^n$, we have
	\begin{equation}
		|A^M_\lambda(f_{K})(x)-f(x)|\leq \frac{L}{4}\left(\frac{\lambda+L/2}{\lambda-L/2}+1\right)r_c^2(x)\,.
	\end{equation}
	Furthermore, in  case $(iii)$, $A^M_\lambda(f_K)$ is an interpolation of $f_K$ in $\mathbb{R}^n$.
\end{itemize}
\end{teo}

\medskip

\begin{nota}
	Theorem \ref{Thm.AprxUF.Cmpct} can be used for the solution of practical problems
	such as salt \& pepper noise removal, in which case K is the compact set given by the part 
	of the image which is noise free.
	On the other hand, Theorem \ref{Thm.AprxCnt} can be applied, for instance, to inpainting 
	of damaged images, where $\Omega$ will be the domain to be inpainted using information about 
	$f_K$, with $K=\Omega^c$. We will discuss these applications of the theory developed here
	in our follow-on paper \cite{ZOC15}.
\end{nota}

\noindent The following corollary of Theorem  \ref{Thm.AprxCnt}  can be thought of as an extension of 
Theorem \ref{Thm.AprxUF.Cmpct}, which concern $A^\infty_\lambda(f_K)$, to the case of finite $M>0$, 
under an extra restriction.

\begin{coro}\label{Cor.AprxCnt}
	Suppose $f:\mathbb{R}^n\mapsto \mathbb{R}$ is bounded and uniformly continuous, with $|f(x)|\leq A_0$
	for some constant $A_0>0$ and all $x\in\mathbb{R}^n$. Assume that
	$f(x)=c_0$ for $|x|\geq r>0$, where $c_0\in \mathbb{R}$ and $r>0$ are constants.
	Let $K\subset \mathbb{R}^n$ be a non-empty compact set satisfying $K\subset \bar{B}(0;\, r)$.
	For $R>r$,  define $K_R:= K\cup B^c(0;\,R)$.
\begin{itemize}
	\item[(i)] Let $\omega$ be the least concave majorant of the modulus of continuity $\omega_f$ of $f$.
		 Assume $a\geq 0$, $b\geq 0$ are such that $\omega(t)\leq at+b$ for $t\geq 0$.
		Then for $\lambda>0$,  $M>A_0+\lambda (R+r)^2$
		and all $x\in \co[K]$, we have
		\begin{equation}
			|A^M_\lambda(f_{K_R})(x)-f(x)|\leq
			\omega\left(r_c(x)+\frac{a}{\lambda}+\sqrt{\frac{2b}{\lambda}}\right)\,.
		\end{equation}
	\item[(ii)] If we further assume that  $f$ is a globally Lipschitz function with Lipschitz constant $L>0$,
		then for $\lambda>0$,  $M>A_0+\lambda (R+r)^2$ and all $x\in \co[K]$, we have
		\begin{equation}
			|A^M_\lambda(f_{K_R})(x)-f(x)|\leq Lr_c(x)+\frac{L^2}{\lambda}\,.
		\end{equation}
	\item[(iii)] If we further assume that  $f$ is a $C^{1,1}$ function such that
		$|Df(x)-Df(y)|\leq L|x-y|$ for all $x,\,y\in \mathbb{R}^n$ and  $L>0$ is a constant,
		then for  $\lambda>L$,
		$M>A_0+\lambda (R+r)^2$ and all $x\in \co[K]$, we have
		\begin{equation}
			|A^M_\lambda(f_{K_R})(x)-f(x)|\leq
				\frac{L}{4}\left(\frac{\lambda+L/2}{\lambda-L/2}+1\right)r_c^2(x)\,.
		\end{equation}
		Furthermore, in  case $(iii)$,  $A^M_\lambda(f_{K_R})$ is an interpolation of $f_K$ in $\R^n$.
\end{itemize}
\end{coro}

%%%%%%%%%%%%%%%%%%%%%%%%%%%%%%%%%%%%%%%%%%%%%%%%%%%%%%%%
\medskip
%%%%%%%%%%%%%%%%%%%%%%%%%%%%%%%%%%%%%%%%%%%%%%%%%%%%%%%%

\begin{nota}
Corollary \ref{Cor.AprxCnt} can be viewed as an extrapolation result for bounded uniformly continuous functions
and for globally Lipschitz functions defined on a compact set. For example, we can define $f$ to be zero outside a large
ball containing $K$ and then apply Theorem \ref{Thm.AprxCnt}. Another reason for such extensions is that if
we simply replace $+\infty$ by a finite $M>0$ in Theorem \ref{Thm.AprxUF.Cmpct}, we are not able
to obtain an error estimate independent of $M$, particularly near the boundary of $\co[K]$.
\end{nota}

%%%%%%%%%%%%%%%%%%%%%%%%%%%%%%%%%%%%%%%%%%%%%%%%%%%%%%%%%%%%%%%%%%%%%%%%%%%%%%%%%%%%%%%%%%%%%%%%%%%%
%%%%%%%%%%%%%%%%%%%%%%%%%%%%%%%%%%%%%%%%%%%%%%%%%%%%%%%%%%%%%%%%%%%%%%%%%%%%%%%%%%%%%%%%%%%%%%%%%%%%
%%%%%%%%%%%%%%%%%%%%%%%%%%%%%%%%%%%%%%%%%%%%%%%%%%%%%%%%%%%%%%%%%%%%%%%%%%%%%%%%%%%%%%%%%%%%%%%%%%%%

\setcounter{equation}{0}
\section{Hausdorff Stability and Regularity}\label{Sec.HauStb}

In this section we establish stability and regularity results for our  approximations. 
The stability properties will be expressed in terms
of a notion of Hausdorff continuity,
 and we first introduce a definition of Hausdorff continuity with respect to closed samples
for transforms of bounded functions on $\mathbb{R}^n$. 

\begin{defi}
\label{def-hauslip}
	Let $\mathcal{B}(\mathbb{R}^n)$ be the class of bounded real-valued functions
	from $\mathbb{R}^n$ to $\mathbb{R}$ and choose a fixed $f  \in \mathcal{B}(\mathbb{R}^n)$ that is uniformly continuous.  	
	A transform $\mathcal{T}:\mathcal{B}(\mathbb{R}^n) \to\mathcal{B}(\mathbb{R}^n)$ is said to be Hausdorff
	continuous with respect to closed sample sets at $f$
	if the mapping $K \mapsto \mathcal{T}(f \chi_K)$ is Hausdorff continuous at each non-empty closed set $K_0 \subset \R^n$, in the sense that for 
	 every $\epsilon>0$, there exists $\delta>0$ such that
	\begin{equation*}
		|\mathcal{T}(f \chi_K)(x)-\mathcal{T}(f \chi_{K_0})(x)|<\epsilon
	\end{equation*}
	for all $x\in\mathbb{R}^n$ whenever $K$ is a non-empty closed set with  $\dist_{\mathcal{H}}(K,K_0)<\delta$, and to be uniformly Hausdorff
	continuous with respect to closed sample sets at $f$ if $\delta$ can be chosen independently of    $K_0$.\\
	\newline A transform $\mathcal{T}:\mathcal{B}(\mathbb{R}^n) \to\mathcal{B}(\mathbb{R}^n)$ is said to be   
		Hausdorff-Lipschitz continuous with respect to closed sample sets at $f$ if the mapping $K \mapsto \mathcal{T}(f\chi_K)$ 
		is Hausdorff-Lipschitz continuous, in the sense that there exists $L>0$ such that
	\begin{equation*}
		|\mathcal{T}(f \chi_K)(x)-\mathcal{T}(f \chi_G)(x)|\leq L\,\dist_{\mathcal{H}}(K,G)
	\end{equation*}
	for all $x\in\mathbb{R}^n$ whenever  $K,\,G\subset\mathbb{R}^n$ are closed sets with $\dist_{\mathcal{H}}(K,G)<\infty$.
\end{defi}

\medskip

\begin{nota}
	It is well known that the Euclidean distance function to a non-empty closed set $K$, i.e. the mapping
	$K\mapsto \dist(\cdot;\,K)$, is  Hausdorff-Lipschitz continuous in the sense that
	$|\dist(x;\,K)-\dist(x;\,G)|\leq \dist_{\mathcal{H}}(K,\,G)$ for all $x\in \R^n$ and non-empty closed sets $K, G$,
	and this is, to our knowledge, the only well-known example of a function satisfying  a Hausdorff-Lipschitz property. 
	A further example, which we will extend here, is given in \cite[Theorem 5.5]{ZOC14a}, where it is shown that
	the mapping $K\mapsto C_{\lambda}^u(\chi_K)$ is  Hausdorff-Lipschitz continuous when $K$ is compact.
\end{nota}

\smallskip

\noindent Our first objective  is to show that the mappings $K\mapsto L^M_\lambda(f_K)$, $K\mapsto U^M_\lambda(f_K)$
and $K\mapsto A^M_\lambda(f_K)$ are uniformly Hausdorff continuous for every bounded uniformly continuous function $f$ with $\sup_{x \in \R^n} |f(x)|<M$. 

\begin{lema}\label{Lem.LU}
	Suppose $f:\mathbb{R}^n\mapsto\mathbb{R}$ is bounded and uniformly continuous, with
	$\sup_{\mathbb{R}^n}|f(x)|\leq A_0$ for some constant $A_0>0$, and
	let $M>A_0$. Then for any fixed $\lambda>0$ and any
	non-empty closed set $K\subset\mathbb{R}^n$,
	\begin{equation*}
		C^l_\lambda(f^M_K)=M-C^u_{\lambda}((M-f)\chi_K)
		\quad\text{and}\quad
		C^u_\lambda(f^{-M}_K)=-M+C^u_{\lambda}((M+f)\chi_K)\,.
	\end{equation*}
\end{lema}
\noindent Now by the assumption that $\sup_{\mathbb{R}^n}|f(x)|\leq A_0<M$, both
$M+f$ and $M-f$ are strictly positive in $\mathbb{R}^n$. Hence, by Definition \ref{Def.LwrUpAv} and Lemma \ref{Lem.LU},
the Hausdorff continuity of the mappings $K\mapsto L^M_\lambda(f_K)$, $K\mapsto U^M_\lambda(f_K)$ and
$K \mapsto A^M_\lambda(f_K)$ 
reduces to the Hausdorff continuity of $K\mapsto C^u_\lambda(g_K)$
for uniformly continuous functions $g:\R^n\mapsto \R$
satisfying 
\begin{equation} \label{bounds-g}
	0<M-A_0\leq g(x)<M\;\;\; \mbox{for all} \;\; x\in \R^n.
\end{equation}

\medskip
\noindent We will thus extend  \cite[Theorem 5.5]{ZOC14a}, that proved  Hausdorff-Lipschitz continuity of 
$K \mapsto C^u_{\lambda}(\chi_K)$ corresponding to the special case $f \equiv 1$ in Definition \ref{def-hauslip}, 
to the general case of any bounded uniformly continuous function $f$.
In the terminology of Definition \ref{def-hauslip}, we will show that the upper transform $C^u_{\lambda}$ is 
uniformly Hausdorff continuous with respect to closed sample sets at
each bounded uniformly continuous function $f$, and is Hausdorff-Lipschitz continuous with respect to closed sample 
sets at each such $f$ that is also globally Lipschitz continuous.

\medskip
\noindent  Motivated by the analysis  in \cite{ZOC14a}, we introduce a squared distance-like function
$D^2_{\lambda,\,f}$,  the upper transform of which is equal to the upper transform $C^u_{\lambda}(f_K)$ of $f_K$ 
and which proves to be a useful tool in the following.

\begin{defi}\label{Def.Dist}
For $f:\mathbb{R}^n\mapsto \mathbb{R}$ with $0<f\leq M$, we define the following distance-like functions 
for a closed set $K\subset \mathbb{R}^n$:
\begin{equation}\label{Eq.Def.SmalDist}
	d_{\lambda,f}(x,\, K)=\inf\left\{|y-x|-\sqrt{\frac{f(y)}{\lambda}},\; y\in K\right\},\quad x\in \mathbb{R}^n;
\end{equation}
and
\begin{equation}\label{Eq.Def.BigDist}
	D_{\lambda,f}(x,\, K)=-\sqrt{\lambda}\min\{ 0,\, d_{\lambda,f}(x,\, K)\}, \quad x\in \mathbb{R}^n.
\end{equation}
\end{defi}

\medskip

\begin{nota}\label{Rmk.Min}
In the definition of $d_{\lambda,f}(x,\, K)$,  if $f$
is  continuous on $\mathbb{R}^n$ and $K$ is closed, the minimum in \eqref{Eq.Def.SmalDist} is attained,
that is, for every $x\in  \mathbb{R}^n$, there exists  $x^K\in K$ such that
$d_{\lambda,f}(x,\, K)=|x^K-x|-\sqrt{f(x^K)/\lambda}$.
Thus if $f$ is continuous, the `$\inf$' in \eqref{Eq.Def.SmalDist} can be replaced by `$\min$'.
\end{nota}

\medskip

\noindent In Theorem \ref{Thm.HauUp} below, we will follow an indirect approach to proving the
Hausdorff continuity of $C_{\lambda}^u(f_K)$ that exploits  the squared function
$D^2_{\lambda,f}(x;\, K)$. Note that it is also possible to give
a direct proof of Hausdorff continuity, avoiding use of $D^2_{\lambda,f}(x;\, K)$, which yields 
a weaker result, namely that for every $\epsilon>0$, there exists $\delta>0$ such that
$|C^u_\lambda(f_K)(x))-C^u_\lambda(f_E)(x))|<\epsilon$
whenever $\dist_{\mathcal{H}}(K,E)<\delta$. Additionally, we can  derive a Hausdorff  continuity 
result  using the Hausdorff continuity of the Moreau envelopes, 
since $C^u_\lambda (f)=M_\lambda(M^\lambda(f)))$ \cite{ZOC14a}, and it can be shown that
\begin{equation*}
	\begin{split}
		&|M_\lambda(f_K)(x)-M_\lambda(f_E)(x)|\leq
			2\lambda\left( \dist_{\mathcal{H}}(K,E)
			+\sqrt{\frac{2M}{\lambda}}\right)\dist_{\mathcal{H}}(K,E)
			+\omega(\dist_{\mathcal{H}}(K,E)),\\[1.5ex]
		&|M^\lambda(f_K)(x)-M^\lambda(f_E)(x)|\leq
		2\lambda\left( \dist_{\mathcal{H}}(K,E)+
		\sqrt{\frac{2M}{\lambda}}\right)\dist_{\mathcal{H}}(K,E)+
		\omega(\dist_{\mathcal{H}}(K,E))\,,
	\end{split}
\end{equation*}
from which a version of Hausdorff continuity of $C_{\lambda}^u(f_K)$ follows.

\bigskip

\noindent Note that the function $D_{\lambda,f}(x;\, K)$ defined in \eqref{Eq.Def.BigDist} 
is a generalisation of $D_\lambda(x;\,K)$  for the characteristic
function $\chi_K$ introduced in \cite[Definition 5.1]{ZOC14a}, since if we take $f\equiv 1$ in \eqref{Eq.Def.BigDist},
we have
\begin{equation*}
	\begin{split}
		D_\lambda(x;\,K)=\max\Big\{0,\, 1-\sqrt{\lambda}\dist(x;\,K)\Big\}
		&=\max\left\{0,\,\sqrt{\lambda}\left(\frac{1}{\sqrt{\lambda}}-\dist(x;\,K)\right)\right\}\\[1.5ex]
		&=\sqrt{\lambda}\max\left\{ 0,\, -\left(\min_{y\in K}|y-x|-\frac{1}{\sqrt{\lambda}}\right)\right\}\\[1.5ex]
		&=-\sqrt{\lambda}\min\Big\{ 0,\, \min_{y\in K}\big\{|y-x|-1/\sqrt{\lambda}\big\}\Big\}=D_{\lambda,f}(x;\,K)\,.
	\end{split}
\end{equation*}
As well as being a tool to investigate the stability of the upper compensated convex transform of characteristic functions,
the geometry-based function $D_{\lambda}(x;\,K)$  has also, for instance,
been used  to find geometric features such as interior corners \cite{ZOC15b}. Hence our
generalised function $D_{\lambda,f}(x,\, K)$ might also have other applications which
 we will  explore elsewhere.

\bigskip
\noindent We start by stating a few preliminary lemmas, the proofs of which are given in Section \ref{Sec.Prfs}.
\medskip

\begin{lema}\label{Lem.Haus.SmlDst}
Suppose $f:\mathbb{R}^n\mapsto \mathbb{R}$ is bounded and uniformly continuous such that
for some constant $M>0$, $0<f(x)\leq M$ for all $x\in \R^n$.
Let $\omega$ be the least concave majorant of the modulus of continuity of $\sqrt{f}$,
which is itself a modulus of continuity.
Let $K,\, E\subset \mathbb{R}^n$ be non-empty  closed
sets with $\dist_{\mathcal{H}}(K,E)<+\infty$. Then for all $x\in \mathbb{R}^n$,
\begin{equation*}
	|d_{\lambda,f}(x,K)-d_{\lambda,f}(x,E)|\leq
	\dist_{\mathcal{H}}(K,E)+\frac{\omega(\dist_{\mathcal{H}}(K,E))}{\sqrt{\lambda}}.
\end{equation*}
\end{lema}

\medskip

\begin{lema}\label{Lem.Haus.Dst}
Under the assumptions of Lemma \ref{Lem.Haus.SmlDst}, we have
\begin{equation}
	|D^2_{\lambda,f}(x,K)-D^2_{\lambda,f}(x,E)|\leq
	2\sqrt{\lambda M}{\rm dist}_{\mathcal{H}}(K,E)+2\sqrt{ M}\omega(\dist_{\mathcal{H}}(K,E)),\quad x\in \mathbb{R}^n.
\end{equation}
\end{lema}

\medskip
\begin{lema}\label{Lem.Haus.UpTr}
	Suppose $\alpha>0$ is a constant and $x_0\in \mathbb{R}^n$, then for $\lambda>0$,
\begin{equation*}
	C^u_\lambda (\alpha\chi_{\{x_0\}})(x)=
	\left\{\begin{array}{ll}
		\displaystyle \lambda(|x-x_0|-\sqrt{\alpha/\lambda})^2,	& \displaystyle |x-x_0|\leq \sqrt{\alpha/\lambda},\\[1.5ex]
		\displaystyle  0,					& \displaystyle |x-x_0|\geq \sqrt{\alpha/\lambda}.
	\end{array}\right.
\end{equation*}
\end{lema}

\medskip
\begin{lema}\label{Lem.Haus.EqUpTr}
Suppose $f$ satisfies the assumptions of Lemma \ref{Lem.Haus.SmlDst}
and $K\subset \mathbb{R}^n$ is closed. Then for $\lambda >0$ and for all $x\in \mathbb{R}^n$,
\begin{equation*}
	C^u_\lambda(f \chi_K)(x)=C^u_\lambda(D^2_{\lambda,f}(\cdot\, ;\,K))(x)\,.
\end{equation*}
\end{lema}

\medskip

\noindent We are now in a position to state our key result on the  Hausdorff stability of the upper compensated convex transform
with respect to closed sample sets at a bounded uniformly continuous positive function $f$.

\medskip
\begin{teo}\label{Thm.HauUp}
Suppose $f:\mathbb{R}^n\mapsto \mathbb{R}$ is bounded and uniformly continuous, with $0<f(x)\leq M$
for some constant $M>0$ and all $x\in \R^n$.
Let $\omega$ be the least concave majorant of the modulus of continuity of $\sqrt{f}$,
which is itself a modulus of continuity.
Let $K,\, E\subset \mathbb{R}^n$ be non-empty  closed
sets with $\dist_{\mathcal{H}}(K,E)<+\infty$.  Then for all $x\in \mathbb{R}^n$,
$C^u_\lambda(f \chi_K)$ is Hausdorff continuous in the sense that 
\begin{equation}
	|C^u_\lambda(f \chi_K)(x)-C^u_\lambda(f \chi_E)(x)|\leq \;
	2\sqrt{\lambda M}{\rm dist}_{\mathcal{H}}(K,E)+2\sqrt{M}
	\, \omega(\dist_{\mathcal{H}}(K,E))\,.
\end{equation}
\end{teo}

\medskip
\begin{coro}\label{Cor.HauCh}
Under the assumptions of Theorem \ref{Thm.HauUp}, if we further assume that
$f:\mathbb{R}^n\mapsto \mathbb{R}$ is a globally Lipschitz continuous function
satisfying $|f(x)-f(y)|\leq L|x-y|$ and $0<\alpha\leq f\leq M$, then
for all $x\in \mathbb{R}^n$, $C^u_\lambda(f \chi_K)$ is Hausdorff-Lipschitz continuous in the sense that 
\begin{equation}
	|C^u_\lambda(f \chi_K)(x)-C^u_\lambda(f \chi_E)(x)|\leq
	\left(2\sqrt{\lambda M}+L\sqrt{\frac{M}{\alpha}}\right)
	{\rm dist}_{\mathcal{H}}(K,E)\,.
\end{equation}
\end{coro}

\bigskip
\noindent We can now easily state the Hausdorff stability theorem for our approximations. 

\begin{teo}\label{Thm.Haus.StbAver}
	Suppose $f:\mathbb{R}^n\mapsto\mathbb{R}$ is bounded and uniformly continuous
	(respectively, globally Lipschitz continuous) and $|f(x)|\leq A_0$ for $x\in \mathbb{R}^n$.
	Then for $M>A_0$, the mappings $K\mapsto L^M_\lambda(f_K)$, $K\mapsto U^M_\lambda(f_K)$,
	$K\mapsto A^M_\lambda(f_K)$ and
	$K\mapsto (SA)^M_\lambda(f_K)$ are uniformly Hausdorff continuous
	(respectively, Hausdorff-Lipschitz continuous).
\end{teo}

\medskip
\noindent We conclude this section by stating  the regularity, or smoothness, of our approximations.
Since our upper, lower and average approximations are globally Lipschitz functions on $\mathbb{R}^n$
and our mixed approximation is a $C^{1,1}$ function, we have the following.

\medskip

\begin{teo}\label{Thm.Reg.Aprx}
	Let  $K \subset \R^n$ be a non-empty closed set and 
	 $f_K:K\subset\mathbb{R}^n\mapsto\mathbb{R}$ be a bounded function  with $|f(x)| <M$ for all $x \in K$. Suppose  $\lambda>0$ and $\tau>0$ are fixed. Then 
\begin{itemize}
	\item[(i)] $U^{M}_\lambda(f_K)$, $L^{M}_\lambda(f_K)$ and $A^{M}_{\lambda}(f_K)$
		are globally Lipschitz continuous on $\mathbb{R}^n$,
		with the Lipschitz constant bounded above by $8\sqrt{M\lambda}$;
	\item[(ii)]  $(SA)^M_{\lambda,\tau}(f_K)$ is a $C^{1,1}$ function on  $\mathbb{R}^n$, and satisfies  
		\begin{equation}\label{Eq.Est.Reg}
			|(SA)^M_{\lambda,\tau}(f_K)(x)-A^{M}_{\lambda}(f_K)(x)|\leq \frac{16M\lambda}{\tau}\;\;\; \mbox{for all}\;\; x \in \R^n.
		\end{equation}
\end{itemize}
\end{teo}

%%%%%%%%%%%%%%%%%%%%%%%%%%%%%%%%%%%%%%%%%%%%%%%%%%%%%%%%%%%%%%%%%%%%%%%%%%%%%%%%%%%%%%%%%%%%%%%%%%%%
%%%%%%%%%%%%%%%%%%%%%%%%%%%%%%%%%%%%%%%%%%%%%%%%%%%%%%%%%%%%%%%%%%%%%%%%%%%%%%%%%%%%%%%%%%%%%%%%%%%%
%%%%%%%%%%%%%%%%%%%%%%%%%%%%%%%%%%%%%%%%%%%%%%%%%%%%%%%%%%%%%%%%%%%%%%%%%%%%%%%%%%%%%%%%%%%%%%%%%%%%

\setcounter{equation}{0}
\section{Proofs of the Main Results}\label{Sec.Prfs}

%%%%%%%%%%%%%%%%%%%%%%%%%%%%%%%%%%%%%%%%%%%%%%%%%%%%%%%%%%%%%%%%%%%%%%%%%%%
%%%%%%%%%%%%%%%%%%%%%%%%%%%%%%%%%%%%%%%%%%%%%%%%%%%%%%%%%%%%%%%%%%%%%%%%%%%
%%%%%%%%%%%%%%%%%%%%%%%%%%%%%%%%%%%%%%%%%%%%%%%%%%%%%%%%%%%%%%%%%%%%%%%%%%%

{\bf Proof of Proposition \ref{Prp.IntSemi}}:  Since $f$ is $2\lambda$-semiconvex, $C^l_\tau(f)(x)=f(x)$
for $x\in \mathbb{R}^n$ and $\tau\geq\lambda$. As $f(x)\leq f^M_K(x)$ for
$x\in \mathbb{R}^n$, we have, for any $x\in K$,
\[
	f(x)=C^l_\tau(f)(x)\leq C^l_\tau(f^M_K)(x)\leq f^M_K(x)=f(x)\,,
\]
where we have applied the ordering property \eqref{Sec2.Eq.Ordering} to show that $C^l_\tau(f)(x)\leq C^l_\tau(f^M_K)(x)$
and \eqref{Sec2.Eq.Ineq} to state that $C^l_\tau(f^M_K)(x)\leq f^M_K(x)$.
Thus $C^l_\tau(f^M_K)$ is an interpolation of $f_K$. Similarly, if $f$ is $2\lambda$-semiconcave,
$C^u_\tau(f^{-M}_K)$ is an interpolation of $f_K$.
\hfill \qed\\

%%%%%%%%%%%%%%%%%%%%%%%%%%%%%%%%%%%%%%%%%%%%%%%%%%%%%%%%%%%%%%%%%%%%%%%%%%%
%%%%%%%%%%%%%%%%%%%%%%%%%%%%%%%%%%%%%%%%%%%%%%%%%%%%%%%%%%%%%%%%%%%%%%%%%%%
%%%%%%%%%%%%%%%%%%%%%%%%%%%%%%%%%%%%%%%%%%%%%%%%%%%%%%%%%%%%%%%%%%%%%%%%%%%

\noindent {\bf Proof of Theorem \ref{Thm.AprxLimBnd}}:  
The first part is immediate from \eqref{Sec3.semicont} and \eqref{Sec3.convtransforms}.
If $f$ is continuous,  it follows from \cite[Theorem 2.3(iii)]{Zha08}  that, uniformly on any compact set,
\[
	\lim_{\lambda\to+\infty} C^l_\lambda(f)(x)=f(x),\qquad
	\lim_{\lambda\to+\infty} C^u_\lambda(f)(x)=f(x)\,,
\]
 whereas
 the ordering properties \eqref{Sec2.Eq.Ordering} and 
\eqref{Sec2.Eq.Ineq} imply that
\[
	C^l_\lambda(f)\leq C^l_\lambda(f^M_K) \leq f^M_K,\quad
	C^u_\lambda(f)\geq C^u_\lambda(f^{-M}_K)\geq f^{-M}_K.
\]
Since $f^M_K = f = f^{-M}_K$ on $K$, it follows that
\[
	\lim_{\lambda\to+\infty} C^l_\lambda(f^M_K)(x)=f(x),\quad
	\lim_{\lambda\to+\infty} C^u_\lambda(f^{-M}_K)(x)=f(x),\quad
	\lim_{\lambda\to+\infty} A^M_\lambda(f_K)(x)=f(x)
\]
uniformly on $ K$, as required. \hfill \qed\\

%%%%%%%%%%%%%%%%%%%%%%%%%%%%%%%%%%%%%%%%%%%%%%%%%%%%%%%%%%%%%%%%%%%%%%%%%%%
%%%%%%%%%%%%%%%%%%%%%%%%%%%%%%%%%%%%%%%%%%%%%%%%%%%%%%%%%%%%%%%%%%%%%%%%%%%
%%%%%%%%%%%%%%%%%%%%%%%%%%%%%%%%%%%%%%%%%%%%%%%%%%%%%%%%%%%%%%%%%%%%%%%%%%%

\noindent {\bf Proof of Proposition \ref{Prp.StrBndAprx}:} By Proposition \ref{Prp.Tra},
without loss of generality, we may assume that $x=0$.
Taking the constant function $\ell(y)=-A_0$ for $y\in\mathbb{R}^n$, we see that
$-A_0=\ell(y)\leq f_K^M(y)+\lambda|y|^2$ so that $-A_0\leq C^l_\lambda(f_K^M)(0)$.

\medskip
\noindent Since $0\in \co[K]$, by Proposition \ref{SecNot.Prp.GlCnvx}, there exist $x_i,\ldots,x_{n+1}\in \mathbb{R}^n$
with $\lambda_i\geq 0$ for $i=1,\ldots, n+1$ such that $\sum^{n+1}_{i=1}\lambda_i=1$
and $\sum^{n+1}_{i=1}\lambda_ix_i=0$. We then have
\[
	C^l_\lambda(f_K^M)(0)=\co[f_K^M+\lambda|\cdot|^2](0)\leq
	\sum^{n+1}_{i=1}\lambda_i(f_K^M(x_i)+\lambda|x_i|^2)
	\leq \sum^{n+1}_{i=1}\lambda_i(A_0+\lambda d^2_K)<M\,.
\]
The proof for the upper transform follows similar arguments. \hfill \qed\\

%%%%%%%%%%%%%%%%%%%%%%%%%%%%%%%%%%%%%%%%%%%%%%%%%%%%%%%%%%%%%%%%%%%%%%%%%%%
%%%%%%%%%%%%%%%%%%%%%%%%%%%%%%%%%%%%%%%%%%%%%%%%%%%%%%%%%%%%%%%%%%%%%%%%%%%
%%%%%%%%%%%%%%%%%%%%%%%%%%%%%%%%%%%%%%%%%%%%%%%%%%%%%%%%%%%%%%%%%%%%%%%%%%%

\noindent {\bf Proof of Theorem \ref{Thm.WMP}:} 
Let $\sup_Kf=A_0$, so that by our assumption \eqref{median}, $\inf_Kf=-A_0$.
Fix $x\in \co[K]$. By Proposition \ref{Prp.Tra}, without loss of generality, we assume that $x=0$.
Notice that $K$ is compact, $C^l_\lambda(f^\infty_K)(0)=C^l_\lambda(\underline{f_K^\infty})(0)$, and
$\underline{f_K^\infty}$ is lower semicontinuous.
Also, $C^u_\lambda(f^{-\infty}_K)(0)=C^u_\lambda(\overline{f_K^{-\infty}})(0)$, and
$\overline{f_K^{-\infty}}$ is upper semicontinuous.
Thus, by Proposition \ref{SecNot.Prp.GlCnvx}, there are two finite generating sets
$K_l=\{x_i^-\}_{i=1}^{k_-}\subset K$ and $K_u=\{x_i^+\}_{i=1}^{k_+}\subset K$, two sets of positive numbers
$\Lambda_l=\{\lambda_i^-\}_{i=1}^{k_-}$ and $\Lambda_u=\{\lambda_i^+\}_{i=1}^{k_+}$
satisfying  $1\leq k_-,\, k_+\leq n+1$, $\sum^{k_-}_{i=1}\lambda_i^-=1$,
$\sum^{k_-}_{i=1}\lambda_i^-x_i^-=0$, $\sum^{k_+}_{i=1}\lambda_i^+=1$,
$\sum^{k_+}_{i=1}\lambda_i^+x_i^+=0$, such that
\[
	\begin{array}{l}
		\displaystyle{\co[\underline{f_K^\infty}+\lambda|\cdot|^2]|(0)
		=\sum^{k_-}_{i=1}\lambda_i^-[\underline{f_K^{\infty}}(x^-_i)+\lambda|x^-_i|^2]}\\
		=
		\displaystyle{\inf\Big\{\sum^{n+1}_{i=1}\lambda_i[\underline{f_K^{\infty}}(x_i)+\lambda|x_i|^2],\;
		x_i\in K_l\cup K_u,\;\lambda_i\geq 0,\; \sum^{n+1}_{i=1}\lambda_i=1,\; \sum^{n+1}_{i=1}\lambda_ix_i=0\Big\}} \\
		\geq \displaystyle{ B_0-A_0},
	\end{array}
\]
where $B_0=\lambda\inf\Big\{\sum^{n+1}_{i=1}\lambda_i|x_i|^2,\; x_i\in K_l\cup K_u,
\;\lambda_i\geq 0,\; \sum^{n+1}_{i=1}\lambda_i=1,\;
\sum^{n+1}_{i=1}\lambda_ix_i=0\Big\}$.
\noindent Likewise
\[
	\co[\underline{f_K^\infty}+\lambda|\cdot|^2](0)\leq B_0+A_0\,,
\]
and thus
\[
	B_0-A_0\leq C^l_\lambda(f^\infty_K)(0)\leq B_0+A_0\,.
\]
On the other hand, we also have, by Proposition \ref{SecNot.Prp.GlCnvx}, that
\[
	\begin{array}{l}
		\displaystyle{
			\co[\lambda|\cdot|^2-\overline{f_K^{-\infty}}](0)
			=\sum^{k_+}_{i=1}\lambda_i^+[\lambda|x^+_i|^2-\overline{f_K}(x^+_i)]
		}\\
		= \displaystyle{
			\inf\left\{\sum^{n+1}_{i=1}\lambda_i[\lambda|x_i|^2-\overline{f_K^{-\infty}}(x_i)],\;
			x_i\in K_l\cup K_u,\;\lambda_i\geq 0,\; \sum^{n+1}_{i=1}\lambda_i=1,\;
			\sum^{n+1}_{i=1}\lambda_ix_i=0\right\}
		} \\
		\geq \displaystyle{B_0-A_0},
\end{array}
\]
and similarly
\[
	\co[\lambda|\cdot|^2-\overline{f_K^{-\infty}}](0)\leq B_0+A_0\,,
\]
so since $C^u_\lambda(f^{-\infty}_K)(0)=-\co[\lambda|\cdot|^2-\overline{f_K^{-\infty}}](0)$,
we obtain
\[
	-B_0-A_0\leq C^u_\lambda(f^{-\infty}_K)(0)\leq -B_0+A_0\,.
\]
Thus
\[
	-A_0\leq A^\infty_\lambda(f_K)(0)\leq A_0\,,
\]
which concludes the proof. \hfill \qed\\

%%%%%%%%%%%%%%%%%%%%%%%%%%%%%%%%%%%%%%%%%%%%%%%%%%%%%%%%%%%%%%%%%%%%%%%%%%%
%%%%%%%%%%%%%%%%%%%%%%%%%%%%%%%%%%%%%%%%%%%%%%%%%%%%%%%%%%%%%%%%%%%%%%%%%%%
%%%%%%%%%%%%%%%%%%%%%%%%%%%%%%%%%%%%%%%%%%%%%%%%%%%%%%%%%%%%%%%%%%%%%%%%%%%

{\bf Proof of Theorem \ref{Thm.AprxUF.Cmpct}.} \textit{Part (i)}:
By Proposition \ref{Prp.Tra}, without loss of generality we  again assume
that $x=0$. Since both $y\mapsto\lambda|y|^2+f^\infty_K(y)$ and $y\mapsto\lambda|y|^2-f^{-\infty}_K(y)$
are coercive and lower semicontinuous, we have, by Proposition \ref{SecNot.Prp.GlCnvx}, that
\[
	C^l_\lambda(f^\infty_K)(0)=\sum^{k_l}_{j=1}\lambda^l_j(\lambda|x_j^l|^2+f(x^l_j)),\quad
	-C^u_\lambda(f^{-\infty}_K)(0)=
	\sum^{k_u}_{j=1}\lambda^u_j(\lambda|x_j^u|^2-f(x^u_j))\,,
\]
where $2\leq k_l,\, k_u\leq n+1$, $\lambda_j^l> 0$, $x_j^l\in K$, $j=1,\dots,k_l$,
$\sum^{k_l}_{j=1}\lambda^l_j=1$, $\sum^{k_l}_{j=1}\lambda^l_jx_j^l=0$; $\lambda_j^u> 0$,
$x_j^u\in K$, $j=1,\dots,k_u$, $\sum^{k_l}_{j=1}\lambda^u_j=1$, $\sum^{k_l}_{j=1}\lambda^u_jx_j^u=0$.

\noindent We also define
\begin{equation*}
\begin{split}
	B_0&=\min\left\{\sum^{n+1}_{k=1}\lambda_k|x_k|^2,\; \lambda_k\geq 0,\; x_k\in K,\;
	k=1,2,\dots,n+1,\sum^{n+1}_{k=1}\lambda_k=1,\;\sum^{n+1}_{k=1}\lambda_kx_k=0\right\}\\[1.5ex]
	&
	=\sum^{m^\ast}_{k=1}\lambda^\ast_k|x^\ast_k|^2\,,
\end{split}
\end{equation*}
for some $2\leq m^\ast\leq n+1$, $\lambda_j^\ast>0$, $x_j^\ast\in K^\ast$ for
$j=1,2,\dots,m^\ast$, $\sum^{m^\ast}_{j=1}\lambda_j^\ast=1$ and
$\sum^{m^\ast}_{j=1}\lambda_j^\ast x_j^\ast=0$,
and let
\begin{equation*}
	\begin{split}
		C_0&=\min\left\{
			\sum^{n+1}_{k=1}\lambda_k|x_k|^2,\;
			\lambda_k\geq 0,\; x_k\in \bar B_{r_c(0)}(0)\cap K,\; k=1,\dots,n+1,\sum^{n+1}_{k=1}\lambda_k=1,\;
			\sum^{n+1}_{k=1}\lambda_kx_k=0
		\right\}
		\\[1.5ex]
		&=\sum^{n+1}_{k=1}\lambda^r_k|x^r_k|^2\,.
	\end{split}
\end{equation*}
Clearly $C_0\leq r^2_c(0)$, and by definition,
\[
	B_0=\sum^{m^\ast}_{k=1}\lambda^\ast_k|x^\ast_k|^2\leq C_0\leq r^2_c(0)\,.
\]
By the Cauchy-Schwarz inequality, we also have
\[
	\sum^{m^\ast}_{k=1}\lambda^\ast_k|x^\ast_k|\leq r_c(0)\,.
\]
Now 
\begin{equation}\label{Eq.AprxUF.Cmpct.1}
	\begin{split}
	 C^l_\lambda(f^\infty_K)(0) & \leq \sum^{m^\ast}_{k=1}\lambda^\ast_k(\lambda|x^\ast_k|^2+f(x^\ast_k))
		 =\lambda B_0+f(0)+\sum^{m^\ast}_{k=1}\lambda^\ast_k(f(x^\ast_k)-f(0))\\
		&\leq \lambda B_0+f(0)+\sum^{m^\ast}_{k=1}\lambda^\ast_k \omega(|x^\ast_k|)
		\leq \lambda B_0+f(0)+\omega\left(\sum^{m^\ast}_{k=1}\lambda^\ast_k|x^\ast_k|\right)\\
		&\leq \lambda B_0+f(0)+\omega(r_c(0)),
	\end{split}
\end{equation}
since $\omega$ is non-decreasing and concave. Furthermore, we also have 
\begin{equation}\label{Eq.AprxUF.Cmpct.2}
	\lambda B_0+f(0)+\omega(r_c(0))\leq
	\lambda r^2_c(0)+f(0)+\omega(r_c(0))\,,
\end{equation}
and 
\begin{equation}\label{Eq.AprxUF.Cmpct.3}
	\begin{split}
		C^l_\lambda(f^\infty_K)(0) &=
					\sum^{k_l}_{j=1}\lambda^l_j(\lambda|x_j^l|^2+f(x^l_j))
					 \geq f(0)+\sum^{k_l}_{j=1}\lambda^l_j(\lambda|x_j^l|^2-|f(x^l_j)-f(0)|)\\
			& \geq f(0)+\sum^{k_l}_{j=1}\lambda^l_j(\lambda|x_j^l|^2-\omega(|x^l_j|)) 
			 \geq f(0)+\sum^{k_l}_{j=1}\lambda^l_j(\lambda|x_j^l|^2-a|x_j^l|-b)\,.
\end{split}
\end{equation}
By comparing \eqref{Eq.AprxUF.Cmpct.1}, \eqref{Eq.AprxUF.Cmpct.2} with \eqref{Eq.AprxUF.Cmpct.3}, it follows that
\begin{equation*}
	f(0)+\sum^{k_l}_{j=1}\lambda^l_j(\lambda|x_j^l|^2-a|x_j^l|-b)\leq \lambda r^2_c(0)+f(0)+\omega(r_c(0)),
\end{equation*}
and hence
\begin{equation*}
	\begin{split}
		\sum^{k_l}_{j=1} \lambda^l_j\left(|x_j^l|-\frac{a}{2\lambda}\right)^2
		    \leq r^2_c(0)+\frac{\omega(r_c(0))}{\lambda}+\frac{a^2}{4\lambda^2}+\frac{b}{\lambda}
		     \leq r^2_c(0)+\frac{a}{\lambda}+\frac{a^2}{4\lambda^2}+\frac{2b}{\lambda} 
		    = \left(r_c(0)+\frac{a}{2\lambda}\right)^2+\frac{2b}{\lambda}\,.
	 \end{split}
\end{equation*}
Here we have used the fact that $\omega(t)\leq at+b$ for $t\geq 0$. Thus by the Cauchy-Schwarz inequality, 
\begin{equation*}
	\sum^{k_l}_{j=1} \lambda^l_j\left||x_j^l|-\frac{a}{2\lambda}\right|\leq \left(\left(r_c(0)+
		 \frac{a}{2\lambda}\right)^2+
		 \frac{2b}{\lambda}\right)^{1/2}\leq
		 r_c(0)+\frac{a}{2\lambda}+\sqrt{2b/\lambda}\,,
\end{equation*}
so that
\begin{equation*}
	\sum^{k_l}_{j=1}\lambda^l_j|x_j^l|\leq r_c(0)+\frac{a}{\lambda}+\sqrt{2b/\lambda}\,.
\end{equation*}
Now 
\begin{equation}\label{Eq.AprxUF.Cmpct.4}
	\begin{split}
		C^l_\lambda(f^\infty_K(0)) &=
			\sum^{k_l}_{j=1}\lambda^l_j(\lambda|x_j^l|^2+f(x^l_j))
			 \geq f(0)+\sum^{k_l}_{j=1}\lambda^l_j(\lambda|x_j^l|^2-|f(x^l_j)-f(0)|)\\
		& \geq f(0)+\lambda B_0-\sum^{k_l}_{j=1}\lambda^l_j\omega(|x^l_j|)
		 \geq f(0)+\lambda B_0-\omega\left(\sum^{k_l}_{j=1}\lambda^l_j|x^l_j|\right)\\
		& \geq f(0)+\lambda B_0 - \omega\left(r_c(0)+a/\lambda+\sqrt{2b/\lambda}\right)\,,
	\end{split}
\end{equation}
and by combining \eqref{Eq.AprxUF.Cmpct.1} and \eqref{Eq.AprxUF.Cmpct.4}, we obtain
\begin{equation*}
	f(0)+\lambda B_0- \omega\left(r_c(0)+a/\lambda+\sqrt{2b/\lambda}\right)\leq
		 C^l_\lambda(f^\infty_K)(0)\leq \lambda B_0+f(0)+\omega(r_c(0))\,.
\end{equation*}
Similarly, we can prove
\begin{equation*}
	f(0)-\lambda B_0-\omega(r_c(0))\leq C^u_\lambda(f^{-\infty}_K)(0)\leq f(0)+\omega\left(r_c(0)+a/\lambda+
	 \sqrt{2b/\lambda}\right)-\lambda B_0\,,
\end{equation*}
  and thus
 \begin{equation*}
	|A^\infty_\lambda(f_K)(0)-f(0)|\leq\frac{1}{2}\left(\omega(r_c(0))+
		 \omega\left(r_c(0)+a/\lambda+\sqrt{2b/\lambda}\right)\right)\,.
\end{equation*}
The proof of Part (i) is thus complete. \\

%%%%%%%%%%%%%%%%%%%%%%%%%%%%%%%%%%%%%%%%%%%%%%%%%%%
\medskip
%%%%%%%%%%%%%%%%%%%%%%%%%%%%%%%%%%%%%%%%%%%%%%%%%%

\noindent \textit{Part (ii)}: We only need to note that in this case, $\omega(t)=Lt$ for $t\geq 0$, taking $a=L$ and $b=0$.
The result then follows. \\

%%%%%%%%%%%%%%%%%%%%%%%%%%%%%%%%%%%%%%%%%%%%%%%%%%%
\medskip
%%%%%%%%%%%%%%%%%%%%%%%%%%%%%%%%%%%%%%%%%%%%%%%%%%

\textit{Part (iii)}:
By Proposition \ref{Prp.Tra} we again assume that $x=0$.
The proof is similar to that of \textit{Part (i)}, and in the following we use
the same notation as in the proof of \textit{Part (i)} for
$\lambda_i^l$, $x_i^l$, $\lambda_j^r$, $x_j^r$ and $\lambda_k^\ast$, $x_k^\ast$.
Thus since
\[
	B_0:=\sum^{m^\ast}_{k=1}\lambda_k^\ast|x_k^\ast|^2\leq \sum^{n+1}_{k=1}\lambda_k^r|x_k^r|^2 \leq r_c^2(0)\,,
\]
we have
\begin{equation}\label{Eq.AprxDif.Cmpct.1}
	\begin{split}
		C^l_\lambda(f^\infty_K)(0)& \leq \sum^{m^\ast}_{k=1}\lambda_k^\ast(f(x_k^\ast)+\lambda |x_k^\ast|^2)
			\; =\; \lambda B_0+f(0)+\sum^{m^\ast}_{k=1}\lambda_k^\ast(f(x_k^\ast)-f(0)-Df(0)\cdot x_k^\ast)\\
			& \leq \lambda B_0+f(0)+\frac{L}{2}\sum^{m^\ast}_{k=1}\lambda_k^\ast|x_k^\ast|^2
			 \;\leq \; \lambda B_0+f(0)+\frac{L}{2}r_c^2(0)\; \leq \; f(0)+\left(\frac{L}{2}+\lambda\right)r_c^2(0)\,,
\end{split}
\end{equation}
and also
\begin{equation}\label{Eq.AprxDif.Cmpct.3}
	\begin{split}
		C^l_\lambda(f^\infty_K)(0)& =
			\sum^{k_l}_{i=1}\lambda_i^l(f(x_i^l)+\lambda |x_i^l|^2)
		 =f(0)+\lambda\sum^{k_l}_{i=1}\lambda_i^l|x_i^l|^2+
			\sum^{k_l}_{i=1}\lambda_i^l(f(x_i^l)-f(0)-Df(0)\cdot x_i^l)\\
		& \geq f(0)+\lambda\sum^{k_l}_{i=1}\lambda_i^l|x_i^l|^2-\frac{L}{2}
			\sum^{k_l}_{i=1}\lambda_i^l|x_i^l|^2
		 =f(0)+\left(\lambda-\frac{L}{2}\right)\sum^{k_l}_{i=1}\lambda_i^l|x_i^l|^2\,.
	\end{split}
\end{equation}
By comparing \eqref{Eq.AprxDif.Cmpct.1} and \eqref{Eq.AprxDif.Cmpct.3}, we then obtain
\begin{equation*}
	f(0)+\left(\lambda-\frac{L}{2}\right)\sum^{k_l}_{i=1}\lambda_i^l|x_i^l|^2
	\leq f(0)+\left(\frac{L}{2}+\lambda\right)r_c^2(0)\,,
\end{equation*}
so that
\begin{equation*}
	\sum^{k_l}_{i=1}\lambda_i^l|x_i^l|^2\leq \frac{\lambda+\frac{L}{2}}{\lambda-\frac{L}{2}}r_c^2(0)\,.
\end{equation*}
Thus from \eqref{Eq.AprxDif.Cmpct.3}, we have
\begin{equation}\label{Eq.AprxDif.Cmpct.4}
	\begin{split}
		C^l_\lambda(f^\infty_K)(0) & = \sum^{k_l}_{i=1}\lambda_i^l(f(x_i^l)+\lambda |x_i^l|^2)\\
			& = f(0)+\lambda\sum^{k_l}_{i=1}\lambda_i^l|x_i^l|^2+
				\sum^{k_l}_{i=1}\lambda_i^l(f(x_i^l)-f(0)-Df(0)\cdot x_i^l)\\
			& \geq f(0)+\lambda B_0-\frac{L}{2}\sum^{k_l}_{i=1}\lambda_i^l|x_i^l|^2
			\; \geq \; f(0)+\lambda B_0-\frac{L}{2}
			\left(\frac{\lambda+\frac{L}{2}}{\lambda-\frac{L}{2}}\right)r_c^2(0)\,.
	\end{split}
\end{equation}
By combining \eqref{Eq.AprxDif.Cmpct.1} and \eqref{Eq.AprxDif.Cmpct.4}, we finally get
\begin{equation*}
	f(0)+\lambda B_0-\frac{L}{2}\left(\frac{\lambda+\frac{L}{2}}{\lambda-\frac{L}{2}}\right)r_c^2(0)
	\leq C^l_\lambda(f^\infty_K)(0)\leq f(0)+\lambda B_0+\frac{L}{2}r_c^2(0)\,.
\end{equation*}
Similarly we can show that
\begin{equation*}
	f(0)-\lambda B_0-\frac{L}{2}r_c^2(0)\leq C^u_\lambda(f^{-\infty}_K)(0)\leq
	f(0)-\lambda B_0+\frac{L}{2}\left(\frac{\lambda+\frac{L}{2}}{\lambda-\frac{L}{2}}\right)r_c^2(0)\,.
\end{equation*}
The conclusion then follows. \hfill \qed\\

%%%%%%%%%%%%%%%%%%%%%%%%%%%%%%%%%%%%%%%%%%%%%%%%%%%
\medskip
%%%%%%%%%%%%%%%%%%%%%%%%%%%%%%%%%%%%%%%%%%%%%%%%%%

\begin{nota}
	From the proof of \textit{Part (i)} of Theorem \ref{Thm.AprxUF.Cmpct} we observe that for a finite $M>0$,
	if $C^l_\lambda(f^M_K)(0)$ can be calculated by using values of $f$ in $K$ ($0\in K$), that is,
	\[
		C^l_\lambda(f^M_K)(0)=\sum^{n+1}_{i=1}\lambda_i\big(f(x_i)+\lambda|x_i|^2\big),
	\]
	with $\lambda_i\geq 0$, $x_i\in K$, $\sum^{n+1}_{i=1}\lambda_i=1$ and
	$\sum^{n+1}_{i=1}\lambda_ix_i=0$, and if a similar result holds for the upper transform,
	then the arguments of the proof of \textit{Part (i)} can go through without any changes.
	However, it is possible that one of the $x_i$'s does not belong to $K$. In this case the
	situation is more complicated.
	In fact, we do not know whether \textit{Part (i)} still holds for a finite $M>0$.
	However, if we extend $f_K$ outside a large ball as zero, we can still derive error bounds
	(see Corollary \ref{Cor.AprxCnt}).
\end{nota}

%%%%%%%%%%%%%%%%%%%%%%%%%%%%%%%%%%%%%%%%%%%%%%%%%%%%%%%%%%%%%%%%%%%%%%%%%%%
%%%%%%%%%%%%%%%%%%%%%%%%%%%%%%%%%%%%%%%%%%%%%%%%%%%%%%%%%%%%%%%%%%%%%%%%%%%
%%%%%%%%%%%%%%%%%%%%%%%%%%%%%%%%%%%%%%%%%%%%%%%%%%%%%%%%%%%%%%%%%%%%%%%%%%%

\noindent {\bf Proof of Theorem \ref{Thm.AprxCnt}.}
\textit{Part (i)}: We first give estimates for $C^l_\lambda(f_K^M)(x)$.
Without loss of generality we assume that $x=0$. Since $y\mapsto f_K^M(y) +\lambda|y|^2$
is lower semicontinuous and coercive, there are $x_i\in\mathbb{R}^n$, $\lambda_i>0$ for
$i=1,2,\dots,m\leq n+1$ such that, $\sum^m_{i=1}\lambda_i=1$, $\sum^m_{i=1}\lambda_ix_i=0$ and
\[
	\co[f_K^M +\lambda|\cdot|^2](0)=\sum^m_{i=1}\lambda_i\big(f_K^M(x_i) +\lambda|x_i|^2\big)\,.
\]
This implies that there is an affine function $\ell(y)$ such that $\ell(y)\leq f_K^M(y) +\lambda|y|^2$
for $y\in\mathbb{R}^n$ and $\ell(x_i)=f_K^M(x_i) +\lambda|x_i|^2$.

\medskip
\noindent We first show that $x_i\in K$ for  $i=1,2,\dots,k$. If this is not the case, there is some
$1\leq i_0\leq k$ such that $x_{i_0}\in\Omega$. Since
$\ell(x_{i_0})=f_K^M(x_{i_0}) +\lambda|x_{i_0}|^2$, $\ell(y)$ is an affine support function of $M+\lambda |y|^2$
at $x_{i_0}$, and hence is the unique tangent plane of
 the function $M+\lambda |y|^2$. Thus $\ell(y)=M+\lambda|x_{i_0}|^2+2\lambda x_{i_0}\cdot y$.

\medskip
\noindent If $x_{i_0}=0$, $\ell(0)=M$, which contradicts the assumption that $M>2A_0+\lambda d^2_\Omega$.
If $x_{i_0}\neq 0$ and $x_{i_0}\in\Omega$, then since $\Omega$ is a bounded domain, there are two
 points $x^\prime_{i_0}$, $x^{\prime\prime}_{i_0}\in\partial\Omega$ and some $0<\alpha<1$, such that
 $x_{i_0}=\alpha x^\prime_{i_0}+(1-\alpha)x^{\prime\prime}_{i_0}$. We also have
 \begin{equation*}
	\begin{split}
		  \alpha \big(f^M_K(x^\prime_{i_0})+\lambda|x^\prime_{i_0}|^2\big)+
		 (1-\alpha)\big(f^M_K(x^{\prime\prime}_{i_0})+\lambda|x^{\prime\prime}_{i_0}|^2\big) & =
		 \alpha \big(f(x^\prime_{i_0})+\lambda|x^\prime_{i_0}|^2\big)+
		 (1-\alpha)\big(f(x^{\prime\prime}_{i_0})+\lambda|x^{\prime\prime}_{i_0}|^2\big)\\
		& \leq A_0+\lambda d_\Omega^2 \\
		& <M\, \leq \, M+\lambda|x_{i_0}|^2		\, = \, f_K^M(x_{i_0}) +\lambda|x_{i_0}|^2\,.
	\end{split}
 \end{equation*}
 Here we have used the fact that $x^\prime_{i_0},\, x^{\prime\prime}_{i_0}\in\partial\Omega$ and
 $0\in\Omega$, so that $|x^\prime_{i_0}-0|\leq d_\Omega$ and $|x^{\prime\prime}_{i_0}-0|\leq d_\Omega$. Thus
 \begin{equation*}
	\begin{split}
		&1=\sum^m_{i=1}\lambda_i=\left(\sum^m_{i=1,i\neq i_0}\lambda_i\right)+\alpha \lambda_{i_0}+(1-\alpha) \lambda_{i_0},\\
		&0=\sum^m_{i=1}\lambda_ix_i=\left(\sum^m_{i=1,i\neq i_0}\lambda_ix_i\right)+ \lambda_{i_0}\alpha x^\prime_{i_0}+
		\lambda_{i_0}(1-\alpha) x^{\prime\prime}_{i_0}\,,
	\end{split}
 \end{equation*}
and
 \begin{eqnarray*}
	 \lefteqn{\co[f_K^M +\lambda|\cdot|^2](0)  =\sum^m_{i=1}\lambda_i[f_K^M(x_i) +\lambda|x_i|^2]}\\
	 		&  > \left(\sum^m_{i=1, i \neq i_0}\lambda_i[f_K^M(x_i) +\lambda|x_i|^2]\right)
			+\alpha \lambda_{i_0} \big( f^M_K(x^\prime_{i_0})+\lambda|x^\prime_{i_0}|^2\big)+
		(1-\alpha) \lambda_{i_0}  \big(f^M_K(x^{\prime\prime}_{i_0})+\lambda|x^{\prime\prime}_{i_0}|^2\big)\,.
 \end{eqnarray*}
But this contradicts the definition of the convex envelope. So $x_i\in K$ for  all $i=1,2,\dots,k$.
The rest of the proof of {\it Part (i)} then follows from a similar argument to that for \textit{Part (i)} of Theorem \ref{Thm.AprxUF.Cmpct}.

\medskip
\noindent For {\it Part (ii)} and {\it Part (iii)}, we can use similar arguments to the proof of \textit{Part (i)}
to show that all $x_i$'s are in $K$, so that the conclusions then follow  from \textit{Part (ii)}
and \textit{Part (iii)} of Theorem \ref{Thm.AprxUF.Cmpct}, respectively.
\hfill \qed\\

%%%%%%%%%%%%%%%%%%%%%%%%%%%%%%%%%%%%%%%%%%%%%%%%%%%%%%%%%%%%%%%%%%%%%%%%%%%
%%%%%%%%%%%%%%%%%%%%%%%%%%%%%%%%%%%%%%%%%%%%%%%%%%%%%%%%%%%%%%%%%%%%%%%%%%%
%%%%%%%%%%%%%%%%%%%%%%%%%%%%%%%%%%%%%%%%%%%%%%%%%%%%%%%%%%%%%%%%%%%%%%%%%%%

\noindent {\bf Proof of Corollary \ref{Cor.AprxCnt}:} The proof  is very similar to that of Theorem \ref{Thm.AprxCnt} and is left to interested readers.
\hfill \qed\\

%%%%%%%%%%%%%%%%%%%%%%%%%%%%%%%%%%%%%%%%%%%%%%%%%%%%%%%%%%%%%%%%%%%%%%%%%%%
%%%%%%%%%%%%%%%%%%%%%%%%%%%%%%%%%%%%%%%%%%%%%%%%%%%%%%%%%%%%%%%%%%%%%%%%%%%
%%%%%%%%%%%%%%%%%%%%%%%%%%%%%%%%%%%%%%%%%%%%%%%%%%%%%%%%%%%%%%%%%%%%%%%%%%%

\noindent {\bf Proof of Lemma \ref{Lem.LU}:}
This lemma  is a direct consequence of the definitions
\eqref{Eq.Def.ExtFnct} of $f^M_K$, $f^{-M}_K$ and the definition
of the upper and lower
compensated convex transforms  \eqref{Eq.Def.UpLwTr}.
\hfill \qed\\

%%%%%%%%%%%%%%%%%%%%%%%%%%%%%%%%%%%%%%%%%%%%%%%%%%%%%%%%%%%%%%%%%%%%%%%%%%%
%%%%%%%%%%%%%%%%%%%%%%%%%%%%%%%%%%%%%%%%%%%%%%%%%%%%%%%%%%%%%%%%%%%%%%%%%%%
%%%%%%%%%%%%%%%%%%%%%%%%%%%%%%%%%%%%%%%%%%%%%%%%%%%%%%%%%%%%%%%%%%%%%%%%%%%

\noindent {\bf Proof of Lemma \ref{Lem.Haus.SmlDst}:}
Fix $x\in\mathbb{R}^n$. For every $\delta> \dist_{\mathcal{H}}(K,E)$, by Remark \ref{Rmk.Min},
there is some $x^E\in E$, such that $d_{\lambda,f}(x,\, E)= |x^E-x|-\sqrt{f(x^E)/\lambda}$.
For  $x^E\in E$, there is some $x^K\in K$ such that $|x^K-x^E|<\delta$. Thus
\begin{equation*}
	\begin{split}
		d_{\lambda,f}(x,\, K)-d_{\lambda,f}(x,\, E)&\leq
			|x^K-x|-\sqrt{f(x^K)/\lambda}-|x^E-x|+\sqrt{f(x^E)/\lambda}\\[1.5ex]
		&\leq |x^K-x^E|+\frac{1}{\sqrt{\lambda}}\omega(|x^K-x^E|)\\[1.5ex]
		&\leq \delta+\frac{\omega(\delta)}{\sqrt{\lambda}}
	\end{split}
\end{equation*}
for all $\delta>\dist_{\mathcal{H}}(K,E)$. Hence,
\[
	d_{\lambda,f}(x,\, K)-d_{\lambda,f}(x,\, E)\leq
	\dist_{\mathcal{H}}(K,E)+\frac{\omega(\dist_{\mathcal{H}}(K,E))}{\sqrt{\lambda}}\,.
\]
Similarly, we can show that
\[
	d_{\lambda,f}(x,\, E)-d_{\lambda,f}(x,\, K)\leq
	\dist_{\mathcal{H}}(K,E)+\frac{\omega(\dist_{\mathcal{H}}(K,E))}{\sqrt{\lambda}}\,,
\]
and conclusion then follows. \hfill \qed\\

%%%%%%%%%%%%%%%%%%%%%%%%%%%%%%%%%%%%%%%%%%%%%%%%%%%%%%%%%%%%%%%%%%%%%%%%%%%
%%%%%%%%%%%%%%%%%%%%%%%%%%%%%%%%%%%%%%%%%%%%%%%%%%%%%%%%%%%%%%%%%%%%%%%%%%%
%%%%%%%%%%%%%%%%%%%%%%%%%%%%%%%%%%%%%%%%%%%%%%%%%%%%%%%%%%%%%%%%%%%%%%%%%%%
\noindent {\bf Proof of Lemma \ref{Lem.Haus.Dst}:}
We have
\begin{equation*}
	|D^2_{\lambda,f}(x,\, K)-D^2_{\lambda,f}(x,\, E)|
	\leq (|D_{\lambda,f}(x,\, K)|+|D_{\lambda,f}(x,\, E)|)|D_{\lambda,f}(x,\, K)-D_{\lambda,f}(x,\, E)|\,.
\end{equation*}
By definition of $D_{\lambda,f}(x,\, K)$, we then have 
if $\min_{y\in K}(|y-x|-\sqrt{f(y)/\lambda}) >0$,
\begin{equation*}
	|D_{\lambda,f}(x,\, K)|=0\,,
\end{equation*}
and if $\min_{y\in K}(|y-x|-\sqrt{f(y)/\lambda})=|x^K-x|- \sqrt{f(x^K)/\lambda}<0$ for some $x^K\in K$, then
\begin{equation*}
	|D_{\lambda,f}(x,\, K)|=\sqrt{\lambda}| |x^K-x|-\sqrt{f(x^K)/\lambda}|=
	\sqrt{\lambda}(\sqrt{f(x^K)/\lambda}-|x^K-x|)\leq \sqrt{f(x^K)} \, \leq \, \sqrt{M}\,.
\end{equation*}
Similarly, we have
\begin{equation*}
	|D_{\lambda,f}(x,\, E)|\leq \sqrt{M}\,.
\end{equation*}

\noindent Next, by the formula $\min\{0,\, a\}=(a- |a|)/2$ for $a\in\mathbb{R}$, we have
\begin{equation*}
	\begin{split}
		|D_{\lambda,f}(x,\, K)-D_{\lambda,f}(x,\, E)|&=\frac{\sqrt{\lambda}}{2}
						\Big|d_{\lambda,f}(x,K)-|d_{\lambda,f}(x,K)|
						-(d_{\lambda,f}(x,E)-|d_{\lambda,f}(x,E)|)\Big|\\[1.5ex]
			& \leq \sqrt{\lambda}\Big|d_{\lambda,f}(x,K)-d_{\lambda,f}(x,E)\Big|\\[1.5ex]
			& \leq \sqrt{\lambda}
				\Big(\dist_{\mathcal{H}}(K,E)+\frac{\omega(\dist_{\mathcal{H}}(K,E))}{\sqrt{\lambda}}\Big)\,.
	\end{split}
\end{equation*}
Thus we obtain
\begin{equation*}
	|D^2_{\lambda,f}(x,\, K)-D^2_{\lambda,f}(x,\, E)|\leq
	2\sqrt{\lambda M}{\rm dist}_{\mathcal{H}}(K,E)+2\sqrt{ M}\omega({\rm dist}_{\mathcal{H}}(K,E))\,,
\end{equation*}
which completes the proof. \hfill \qed\\

%%%%%%%%%%%%%%%%%%%%%%%%%%%%%%%%%%%%%%%%%%%%%%%%%%%%%%%%%%%%%%%%%%%%%%%%%%%
%%%%%%%%%%%%%%%%%%%%%%%%%%%%%%%%%%%%%%%%%%%%%%%%%%%%%%%%%%%%%%%%%%%%%%%%%%%
%%%%%%%%%%%%%%%%%%%%%%%%%%%%%%%%%%%%%%%%%%%%%%%%%%%%%%%%%%%%%%%%%%%%%%%%%%%
\noindent {\bf Proof of Lemma \ref{Lem.Haus.UpTr}:}
The proof of this lemma is an easy exercise and is omitted here. \hfill \qed\\

%%%%%%%%%%%%%%%%%%%%%%%%%%%%%%%%%%%%%%%%%%%%%%%%%%%%%%%%%%%%%%%%%%%%%%%%%%%
%%%%%%%%%%%%%%%%%%%%%%%%%%%%%%%%%%%%%%%%%%%%%%%%%%%%%%%%%%%%%%%%%%%%%%%%%%%
%%%%%%%%%%%%%%%%%%%%%%%%%%%%%%%%%%%%%%%%%%%%%%%%%%%%%%%%%%%%%%%%%%%%%%%%%%%

\noindent {\bf Proof of Lemma \ref{Lem.Haus.EqUpTr}:}
We first show that
\begin{equation}\label{Eq.Lem.Haus.EqUpTr.01}
	f(x)\chi_K(x)\leq D^2_{\lambda,f}(x,\, K)
\end{equation}
for all $x\in \R^n$, so that by \eqref{Sec2.Eq.Ineq},
\begin{equation}\label{Eq.Lem.Haus.EqUpTr.02}
	C^u_\lambda(f\chi_K)(x)\leq C^u_\lambda(D^2_{\lambda,f}(\cdot,\, K))(x)
\end{equation}
 for all $x\in\mathbb{R}^n$.
If $x\notin K$, clearly, $\chi_K(x)f(x)=0\leq D^2_{\lambda,f}(x,\, K)$.
If $x\in K$, since
\[
	d_{\lambda,f}(x,K)=\min_{y\in K}(|y-x|-\sqrt{f(y)/\lambda})\leq -\sqrt{f(x)/\lambda}\,<\, 0\,,
\]
we have
\[
	D_{\lambda,f}(x,K)=- \sqrt{\lambda}\min\{ 0,\, d_{\lambda,f}(x,K)\}=-\sqrt{\lambda}\,d_{\lambda,f}(x,K)\, \geq \sqrt{\lambda}
	 \sqrt{f(x)/\lambda}=\sqrt{f(x)}\,,
\]
 and thus $D^2_{\lambda,f}(x,K)\geq f(x)$. Therefore \eqref{Eq.Lem.Haus.EqUpTr.01} holds for all $x \in \R^n$, from which \eqref{Eq.Lem.Haus.EqUpTr.02} follows.

%\noindent Since \eqref{Eq.Lem.Haus.EqUpTr.01} is satisfied, clearly, we have
%$C^u_\lambda(f\chi_K)(x)\leq C^u_\lambda(D^2_{\lambda,f}(\cdot,K))(x)$ for all $x\in\mathbb{R}^n$.

\bigskip
\noindent Next we show that the opposite inequality, $C^u_\lambda(D^2_{\lambda,f}(\cdot,K))(x)\leq C^u_\lambda(f\chi_K)(x)$, also holds. If
\begin{equation} \label{Eq.Lem.Haus.EqUpTr.03}
	d_{\lambda,f}(x,K)=\min_{y\in K}(|y-x|-\sqrt{f(y)/\lambda})>0,
\end{equation}
then by definition, $D_{\lambda,f}(x,K)=0$, and hence $D^2_{\lambda,f}(x,K)=0$.
We show in this case that
\[
	C^u_\lambda(D^2_{\lambda,f}(\cdot,K))(x)=0\leq C^u_\lambda(f\chi_K)(x)\,.
\]
We will consider the function $z\mapsto \lambda|z-x|^2-D^2_{\lambda,f}(z,K)$ for $z\in \mathbb{R}^n$
and show that the value of the convex envelope of this function at $x$ is zero.
Consider the affine function $\ell(z)=0$ for all $z\in \mathbb{R}^n$ and show that
\begin{equation}\label{Eq.Lem.Haus.EqUpTr.04}
	0=\ell(x)=\Big(\lambda|z-x|^2-D^2_{\lambda,f}(z,K)\Big)|_{z=x},
 \end{equation}
 and
 \begin{equation}\label{Eq.Lem.Haus.EqUpTr.05}
	0=\ell(z)\leq\lambda|z-x|^2-D^2_{\lambda,f}(z,K),\quad z\in \mathbb{R}^n.
\end{equation}

\medskip
\noindent Equality \eqref{Eq.Lem.Haus.EqUpTr.04} is obvious as $[\lambda|z-x|^2-D^2_{\lambda,f}(z,K)]|_{z=x}=-D^2_{\lambda,f}(x,K)=0$. Now we prove \eqref{Eq.Lem.Haus.EqUpTr.05}, that is,
$0\leq \lambda|z-x|^2-D^2_{\lambda,f}(z,K)$,  which is equivalent to
\begin{equation}\label{Eq.Lem.Haus.EqUpTr.06}
	D^2_{\lambda,f}(z,K)\leq \lambda|z-x|^2,\quad z\in \mathbb{R}^n.
\end{equation}
If $d_{\lambda,f}(z,K)\geq 0$, then $D^2_{\lambda,f}(z,K)=0$, hence \eqref{Eq.Lem.Haus.EqUpTr.06} holds. 
If $d_{\lambda,f}(z,K)<0 $, then $D^2_{\lambda,f}(z,K)=\lambda d^2_{\lambda,f}(z,K)$. We need to show that
$\lambda(\min_{y\in K}(|y-z|-\sqrt{f(y)/\lambda})^2\leq \lambda|z-x|^2,$
which is equivalent to
$-\min_{y\in K}(|y-z|-\sqrt{f(y)/\lambda})\leq |z-x|$, which is in turn equivalent to
\begin{equation}\label{Eq.Lem.Haus.EqUpTr.07}
	|z-x|+\min_{y\in K}(|y-z|-\sqrt{f(y)/\lambda})\geq 0.
\end{equation}
By the triangle inequality and \eqref{Eq.Lem.Haus.EqUpTr.03}, we have
\[
	\begin{split}
		|z-x|+\min_{y\in K}(|y-z|-\sqrt{f(y)/\lambda})& =\min_{y\in K}(|z-x|+|y-z|-\sqrt{f(y)/\lambda})\\[1.5ex]
				& \geq \min_{y\in K}(|y-x|-\sqrt{f(y)/\lambda})=d_{\lambda,f}(x,K)>0.
	  \end{split}
\]
Thus \eqref{Eq.Lem.Haus.EqUpTr.05} holds. Therefore
\begin{equation*}
	0=\co[\lambda|\cdot-x|^2-D^2_{\lambda,f}(\cdot,K)](x) =-C^u_\lambda(D^2_{\lambda,f}(\cdot,K))(x),
\end{equation*}
which implies
\begin{equation}\label{Eq.Lem.Haus.EqUpTr.08}
	D^2_{\lambda,f}(x,K)\leq C^u_\lambda(D^2_{\lambda,f}(\cdot,K))(x)=0\leq C^u_\lambda(f\chi_K)(x)\,.
\end{equation}
Finally, we consider the case
\begin{equation}\label{Eq.Lem.Haus.EqUpTr.09}
	 d_{\lambda,f}(x,K)=\min_{y\in K}(|y-x|-\sqrt{f(y)/\lambda})=|x^K-x|-\sqrt{f(x^K)/\lambda} \, < \, 0,
\end{equation}
where $x^K\in K$ is the minimum point. Now we consider the function
$f(y)\chi_{\{x^K\}}(y)$ for $y\in \mathbb{R}^n$. By Lemma \ref{Lem.Haus.UpTr}, we have
\begin{equation*}
	C^u_\lambda (f\chi_{\{x^K\}})(y)=
		\left\{\begin{array}{ll}
	\displaystyle		\lambda\left(|y-x^K|-\sqrt{f(x^K)/\lambda}\right)^2,&
					\displaystyle |y-x^K|\leq \sqrt{f(x^K)/\lambda},\\
	\displaystyle		0,& \displaystyle |y-x^K|\geq \sqrt{f(x^K)/\lambda}.
			\end{array}\right.
\end{equation*}
In particular, since $x^K\in K$, we have  $f(y)\chi_{\{x^K\}}(y)\leq f(y)\chi_K(y)$ for all $y\in\mathbb{R}^n$, so that
by \eqref{Sec2.Eq.Ineq},
\begin{equation*}
	C^u_\lambda (f\chi_{\{x^K\}})(y)\leq C^u_\lambda (f\chi_K)(y)\quad\text{for all }y\in\mathbb{R}^n\,.
\end{equation*}	
By our assumption \eqref{Eq.Lem.Haus.EqUpTr.09}, we also have $|x^K-x|<\sqrt{f(x^K)/\lambda}$, thus
\begin{equation*}
	C^u_\lambda (f\chi_{\{x^K\}})(x)= \lambda\left(|x-x^K|-\sqrt{f(x^K)/\lambda}\right)^2=D^2_{\lambda,f}(x,K)
\end{equation*}
as $d_{\lambda,f}(x,K)<0$. Thus, in this case, 
$D^2_{\lambda,f}(x,K)=C^u_\lambda (f\chi_{\{x^K\}})(x)\leq C^u_\lambda (f\chi_K)(x)$.
By combining this case and \eqref{Eq.Lem.Haus.EqUpTr.08}, we have, for all $x\in\mathbb{R}^n$, that
$D^2_{\lambda,f}(x,K)\leq C^u_\lambda (f\chi_K)(x)$, so that
\begin{equation*}
	C^u_\lambda (D^2_{\lambda,f}(\cdot,K))(x)\leq
	C^u_\lambda (C^u_\lambda (f\chi_K))(x)=C^u_\lambda (f\chi_K)(x)\,.
\end{equation*}
Since  the opposite inequality \eqref{Eq.Lem.Haus.EqUpTr.02} also holds,
we have 
\begin{equation*}
	C^u_\lambda (D^2_{\lambda,f}(\cdot,K))(x)= C^u_\lambda (f\chi_K)(x)
\end{equation*}
for all $x\in \mathbb{R}^n$, which completes the proof. \hfill \qed\\

%%%%%%%%%%%%%%%%%%%%%%%%%%%%%%%%%%%%%%%%%%%%%%%%%%%%%%%%%%%%%%%%%%%%%%%%%%%
%%%%%%%%%%%%%%%%%%%%%%%%%%%%%%%%%%%%%%%%%%%%%%%%%%%%%%%%%%%%%%%%%%%%%%%%%%%
%%%%%%%%%%%%%%%%%%%%%%%%%%%%%%%%%%%%%%%%%%%%%%%%%%%%%%%%%%%%%%%%%%%%%%%%%%%

\noindent {\bf Proof of Theorem \ref{Thm.HauUp}:} By Lemma \ref{Lem.Haus.EqUpTr}, we only need to prove
 \begin{equation}\label{Eq.HauUp}
	 |C^u_\lambda (D^2_{\lambda,f}(\cdot,K))(x)-C^u_\lambda (D^2_{\lambda,f}(\cdot,E))(x)|\leq
	 2\sqrt{\lambda M}{\rm dist}_{\mathcal{H}}(K,E)+2\sqrt{ M}\omega( \dist_{\mathcal{H}}(K,E))\,.
\end{equation}
By Lemma \ref{Lem.Haus.Dst} we have, for all $x\in\mathbb{R}^n$ that
\[
	|D^2_{\lambda,f}(x,K)-D^2_{\lambda,f}(x,E)|\leq
	2\sqrt{\lambda M}{\rm dist}_{\mathcal{H}}(K,E)+2\sqrt{ M}\omega( \dist_{\mathcal{H}}(K,E))\,.
\]
 Thus
 \[
	\begin{split}
		& D^2_{\lambda,f}(x,E)-2\sqrt{\lambda M}{\rm dist}_{\mathcal{H}}(K,E)-2\sqrt{ M}\omega( \dist_{\mathcal{H}}(K,E))
		\leq D^2_{\lambda,f}(x,K)\\[1.5ex]
		&	\leq D^2_{\lambda,f}(x,E)+2\sqrt{\lambda M}{\rm dist}_{\mathcal{H}}(K,E)+2\sqrt{ M}\omega( \dist_{\mathcal{H}}(K,E))
	\end{split}
\]
 for all $x\in \mathbb{R}^n$. By the ordering and the affine covariance properties of compensated convex transforms, we have
\[
	\begin{split}
		& C^u_\lambda (D^2_{\lambda,f}(\cdot,E))(x)-2\sqrt{\lambda M}{\rm dist}_{\mathcal{H}}(K,E)-
			2\sqrt{ M}\omega( \dist_{\mathcal{H}}(K,E))
			\leq C^u_\lambda (D^2_{\lambda,f}(\cdot,K))(x)\\[1.5ex]
		& \leq
		C^u_\lambda (D^2_{\lambda,f}(\cdot,E))(x)+
		2\sqrt{\lambda M}{\rm dist}_{\mathcal{H}}(K,E)+2\sqrt{ M}\omega( \dist_{\mathcal{H}}(K,E)).
	 \end{split}
 \]
 Hence \eqref{Eq.HauUp} follows.  \hfill \qed\\

%%%%%%%%%%%%%%%%%%%%%%%%%%%%%%%%%%%%%%%%%%%%%%%%%%%%%%%%%%%%%%%%%%%%%%%%%%%
%%%%%%%%%%%%%%%%%%%%%%%%%%%%%%%%%%%%%%%%%%%%%%%%%%%%%%%%%%%%%%%%%%%%%%%%%%%
%%%%%%%%%%%%%%%%%%%%%%%%%%%%%%%%%%%%%%%%%%%%%%%%%%%%%%%%%%%%%%%%%%%%%%%%%%%

\noindent {\bf Proof of Corollary \ref{Cor.HauCh}:}
 This is a direct consequence of Theorem \ref{Thm.HauUp} where
 $\omega(t) =Lt/(2\sqrt{\alpha})$ for $t>0$, with $L\geq 0$ the Lipschitz constant of $f$, 
 since  $|\sqrt{f(x)}-\sqrt{f(y)}|\leq L|x-y|/(2\sqrt{\alpha})$.  \hfill \qed\\

%%%%%%%%%%%%%%%%%%%%%%%%%%%%%%%%%%%%%%%%%%%%%%%%%%%%%%%%%%%%%%%%%%%%%%%%%%%
%%%%%%%%%%%%%%%%%%%%%%%%%%%%%%%%%%%%%%%%%%%%%%%%%%%%%%%%%%%%%%%%%%%%%%%%%%%
%%%%%%%%%%%%%%%%%%%%%%%%%%%%%%%%%%%%%%%%%%%%%%%%%%%%%%%%%%%%%%%%%%%%%%%%%%%

\noindent {\bf Proof of Theorem \ref{Thm.Haus.StbAver}:}
The cases of $L^M_\lambda(f_K)$,
$U^M_\lambda(f_K)$ and  $A^M_\lambda(f_K)$ are direct consequences
Definition \ref{Def.LwrUpAv}, Lemma \ref{Lem.LU}, Theorem \ref{Thm.HauUp} and Corollary \ref{Cor.HauCh}.
Note that for the Hausdorff-Lipschitz continuity properties, the assumption that $M>A_0$ ensures that 
the uniform positivity assumption in Corollary \ref{Cor.HauCh} is satisfied by both  $M-f$ and $M+f$.
For the mixed average approximation
$(SA)^M_{\lambda,\tau}(f_K)=\frac{1}{2}(C^u_\tau(C^l_\lambda(f^M_K))+C^l_\tau(C^u_\lambda(f^{-M}_K)))$,  
we use  \eqref{Sec2.Eq.Ineq}.
Since
$|C^l_\lambda(f^M_G)(x)-C^l_\lambda(f^M_K)(x)|<\epsilon$ for all $x\in\mathbb{R}^n$
with $\epsilon=2\sqrt{\lambda M}\dist_{\mathcal{H}}(G,K)+2\sqrt{M}\omega\left(\dist_{\mathcal{H}}(G,K)\right)$,
we have
\[
	C^l_\lambda(f^M_K)(x)-\epsilon<C^l_\lambda(f^M_G)(x)<C^l_\lambda(f^M_K)(x)+\epsilon\,,
\]
and hence
\[
	|C^u_\tau(C^l_\lambda(f^M_G))(x)-C^u_\tau(C^l_\lambda(f^M_K))(x)| \,<\,  \epsilon\,,
\]
since $C^u_{\tau}(C^l_\lambda(f^M_K) \pm \epsilon) = C^u_{\tau}(C^l_\lambda(f^M_K)) \pm \epsilon$.
Similarly,
\[
	|C^l_\tau(C^u_\lambda(f^{-M}_G))(x)-C^l_\tau(C^u_\lambda(f^{-M}_K))(x)|  \,<\,  \epsilon
\]
since $|C^u_\lambda(f^{-M}_G)(x)-C^u_\lambda(f^{-M}_K)(x)|<\epsilon$.
The proof for $(SA)^M_{\lambda,\tau}(f_K)$ then follows.
The proof for the Lipschitz case is similar, using arguments from Lemma \ref{Lem.LU}
and Corollary \ref{Cor.HauCh}.
\hfill \qed\\

%%%%%%%%%%%%%%%%%%%%%%%%%%%%%%%%%%%%%%%%%%%%%%%%%%%%%%%%%%%%%%%%%%%%%%%%%%%%%%%%%%%%%%%%%%%%%%
%%%%%%%%%%%%%%%%%%%%%%%%%%%%%%%%%%%%%%%%%%%%%%%%%%%%%%%%%%%%%%%%%%%%%%%%%%%%%%%%%%%%%%%%%%%%%%
%%%%%%%%%%%%%%%%%%%%%%%%%%%%%%%%%%%%%%%%%%%%%%%%%%%%%%%%%%%%%%%%%%%%%%%%%%%%%%%%%%%%%%%%%%%%%%

\noindent {\bf Proof of Theorem \ref{Thm.Reg.Aprx}:}
Part (i) and the error estimate \eqref{Eq.Est.Reg}
follow from \cite[Theorem 3.13]{ZOC14a}.
The fact that mixed transforms are $C^{1,1}$ is a consequence of  \cite[Theorem 2.1(iv), Theorem 4.1(ii)]{Zha08}.
Note that this latter regularity result also follows from the fact  that if $g$ is both $2\lambda$-semiconvex and $2\lambda$-semiconcave,
then $g$ is a $C^{1,1}$ function \cite[Corollary 3.3.8]{CS04}.
\hfill \qed\\

%%%%%%%%%%%%%%%%%%%%%%%%%%%%%%%%%%%%%%%%%%%%%%%%%%%%%%%%%%%%%%%%%%%%%%%%%%%%%%%%%%%%%%%%%%%%%%
%%%%%%%%%%%%%%%%%%%%%%%%%%%%%%%%%%%%%%%%%%%%%%%%%%%%%%%%%%%%%%%%%%%%%%%%%%%%%%%%%%%%%%%%%%%%%%
%%%%%%%%%%%%%%%%%%%%%%%%%%%%%%%%%%%%%%%%%%%%%%%%%%%%%%%%%%%%%%%%%%%%%%%%%%%%%%%%%%%%%%%%%%%%%%

\noindent {\bf Acknowledgements:}
The authors are grateful to an anonymous referee, for pointing out 
the work on the proximal average. KZ wishes to thank The University of Nottingham for its support,
EC is grateful for the financial support of the College of Science, Swansea University,
and
AO acknowledges the financial support of the Argentinean Agency through the
Project Prestamo BID PICT PRH 30 No 94 and the National University of Tucum\'{a}n through the project PIUNT E527.

%%%%%%%%%%%%%%%%%%%%%%%%%%%%%%%%%%%%%%%%%%%%%%%%%%%%%%%%%%%%%%%%%%%%%%%%%%%%%%%%%%%%%%%%%%%%%%
%%%%%%%%%%%%%%%%%%%%%%%%%%%%%%%%%%%%%%%%%%%%%%%%%%%%%%%%%%%%%%%%%%%%%%%%%%%%%%%%%%%%%%%%%%%%%%
%%%%%%%%%%%%%%%%%%%%%%%%%%%%%%%%%%%%%%%%%%%%%%%%%%%%%%%%%%%%%%%%%%%%%%%%%%%%%%%%%%%%%%%%%%%%%%

\end{document}